\documentclass{article}

\usepackage[titletoc,title]{appendix}
\usepackage[cp1250]{inputenc}
\usepackage[IL2]{fontenc}
\usepackage{yfonts,fancyhdr}
\usepackage{a4wide}
\usepackage[english]{babel}
\usepackage{euscript}
\usepackage{amstext,amsbsy,amscd,amssymb}
\usepackage{amsmath,enumerate}
\usepackage{amsfonts}
\usepackage{mathrsfs,tensor,xy}
\usepackage{graphics}
\usepackage{microtype}
\include{diagxy}
\allowdisplaybreaks

 \input xypic
 \input epsf
 \input xy
 \xyoption{all}

\let\rarr=\rightarrow

\let\lrarr=\longrightarrow
\let\veps=\varepsilon
\let\mcal=\mathcal
\let\mfrak=\mathfrak
\let\eus=\EuScript

\def\N{\mathbb{N}}
\def\Z{\mathbb{Z}}
\def\R{\mathbb{R}}
\def\C{\mathbb{C}}
\def\imag{{\rm i}}

\def\Sp{\mathop {\rm Sp} \nolimits}

\def\End{\mathop {\rm End} \nolimits}
\def\Hom{\mathop {\rm Hom} \nolimits}

\def\Ind{\mathop {\rm Ind} \nolimits}

\def\ad{\mathop {\rm ad} \nolimits}

\def\GL{\mathop {\rm GL} \nolimits}

\def\Sp{\mathop {\rm Sp} \nolimits}
\def\Mp{\mathop{\rm Mp}\nolimits}

\def\diag{\mathop {\rm diag} \nolimits}
\def\id{\mathop {\rm id} \nolimits}

\def\Ad{\mathop {\rm Ad} \nolimits}

\def\Sol{\mathop {\rm Sol} \nolimits}

\def\tr{\mathop {\rm tr} \nolimits}
\def\htt{\mathop {\rm ht} \nolimits}

\newcommand{\mP}{\mathbb{P}}
\newcommand{\mS}{\mathbb{S}}
\newcommand{\mV}{\mathbb{V}}
\newcommand{\mC}{\mathbb{C}}

\newcommand{\mN}{\mathbb{N}}
\newcommand{\lC}{\mathcal{C}}
\newcommand{\mZ}{\mathbb{Z}}
\newcommand{\mR}{\mathbb{R}}
\newcommand{\mW}{\mathbb{W}}

\newcommand{\gog}{\mathfrak{g}}
\newcommand{\gop}{\mathfrak{p}}

\newcommand{\gol}{\mathfrak{l}}

\newsymbol\squares 1003

\long\def\proof #1{\noindent \emph{Proof.}\ #1 \hfill $\squares$
\medskip}

\newcounter{num}[section]
\numberwithin{equation}{section}
\numberwithin{num}{section}

\long\def\definition #1 {\refstepcounter{num} \noindent {\bf
Definition \thenum.} #1

\medskip}

\long\def\theorem #1{\refstepcounter{num} \noindent {\bf Theorem
\thenum.} #1

\medskip}

\long\def\lemma #1{\refstepcounter{num}  \noindent {\bf Lemma
\thenum.} #1

\medskip}

\long\def\proposition #1{\refstepcounter{num}  \noindent {\bf
Proposition \thenum.} #1

\medskip}

\long\def\corollary #1{\refstepcounter{num}  \noindent {\bf
Corollary \thenum.} #1

\medskip}

\newcommand*\riso{%
  \xrightarrow[]{\raisebox{-0.25em}{\smash{\ensuremath{\sim}}}}%
}

\makeatletter
\newcommand*\if@single[3]{%
  \setbox0\hbox{${\mathaccent"0362{#1}}^H$}%
  \setbox2\hbox{${\mathaccent"0362{\kern0pt#1}}^H$}%
  \ifdim\ht0=\ht2 #3\else #2\fi
  }
\newcommand*\rel@kern[1]{\kern#1\dimexpr\macc@kerna}
\newcommand*\widebar[1]{\@ifnextchar^{{\wide@bar{#1}{0}}}{\wide@bar{#1}{1}}}
\newcommand*\wide@bar[2]{\if@single{#1}{\wide@bar@{#1}{#2}{1}}{\wide@bar@{#1}{#2}{2}}}
\newcommand*\wide@bar@[3]{%
  \begingroup
  \def\mathaccent##1##2{%
    \if#32 \let\macc@nucleus\first@char \fi
    \setbox\z@\hbox{$\macc@style{\macc@nucleus}_{}$}%
    \setbox\tw@\hbox{$\macc@style{\macc@nucleus}{}_{}$}%
    \dimen@\wd\tw@
    \advance\dimen@-\wd\z@
    \divide\dimen@ 3
    \@tempdima\wd\tw@
    \advance\@tempdima-\scriptspace
    \divide\@tempdima 10
    \advance\dimen@-\@tempdima
    \ifdim\dimen@>\z@ \dimen@0pt\fi
    \rel@kern{0.6}\kern-\dimen@
    \if#31
      \overline{\rel@kern{-0.6}\kern\dimen@\macc@nucleus\rel@kern{0.4}\kern\dimen@}%
      \advance\dimen@0.4\dimexpr\macc@kerna
      \let\final@kern#2%
      \ifdim\dimen@<\z@ \let\final@kern1\fi
      \if\final@kern1 \kern-\dimen@\fi
    \else
      \overline{\rel@kern{-0.6}\kern\dimen@#1}%
    \fi
  }%
  \macc@depth\@ne
  \let\math@bgroup\@empty \let\math@egroup\macc@set@skewchar
  \mathsurround\z@ \frozen@everymath{\mathgroup\macc@group\relax}%
  \macc@set@skewchar\relax
  \let\mathaccentV\macc@nested@a
  \if#31
    \macc@nested@a\relax111{#1}%
  \else
    \def\gobble@till@marker##1\endmarker{}%
    \futurelet\first@char\gobble@till@marker#1\endmarker
    \ifcat\noexpand\first@char A\else
      \def\first@char{}%
    \fi
    \macc@nested@a\relax111{\first@char}%
  \fi
  \endgroup
}
\makeatother

\newcommand\rsmraise[1]{%
  \ifx#1\displaystyle .8\else
    \ifx#1\textstyle .8\else
      \ifx#1\scriptstyle .6\else
        .45%
      \fi
    \fi
  \fi}

\title{Differential invariants on symplectic spinors in contact projective geometry}

\author{Libor Křižka, Petr Somberg}

\AtEndDocument{\bigskip{\footnotesize%
  (L.\,Křižka) \textsc{Charles University, Faculty of Mathematics and Physics,
	Mathematical Institute, \\ Sokolovská 83, 180\,00 Praha 8, Czech Republic} \par
  \textit{E-mail address}: \texttt{krizka.libor@gmail.com} \par
  \addvspace{\medskipamount}
  (P.\,Somberg) \textsc{Charles University, Faculty of Mathematics and Physics,
	Mathematical Institute, \\ Sokolovská 83, 180\,00 Praha 8, Czech Republic} \par
  \textit{E-mail address}: \texttt{somberg@karlin.mff.cuni.cz} \par
}}

\begin{document}
\date{}
\maketitle

\begin{abstract} We present a complete classification and the construction of
$\Mp(2n+2,\R)$-equivariant differential operators acting on the principal series
representations, associated to the contact projective geometry on $\mathbb{RP}^{2n+1}$
and induced from
the irreducible $\Mp(2n,\R)$-submodules of the Segal--Shale--Weil
representation twisted by a one-parameter family of characters.
The proof is based on the classification of homomorphisms of generalized Verma modules for the
Segal--Shale--Weil representation twisted by a one-parameter family of characters, together with
a generalization of the well-known duality between homomorphisms of generalized Verma modules
and equivariant differential operators in the category of inducing smooth admissible modules.

\medskip

\medskip

\end{abstract}

\thispagestyle{empty}

\tableofcontents


\section*{Introduction}
\addcontentsline{toc}{section}{Introduction}

The idea of a symplectic analogue of spinor fields in Riemannian
(or conformal) geometry is based on the use of the Segal--Shale--Weil representation and originated
in the work by B.\,Kostant, \cite{Kostant1974}. The
geometrical and functional theoretical properties
of the associated first order differential invariant,
the symplectic Dirac operator $D_s$, were consequently studied on
symplectic manifolds, cf.\ \cite{Habermann-Habermann2006} and the references
therein. However, there is another useful realization of symplectic differential
invariants like $D_s$, which is based on an extension of their
symplectic symmetry group. In particular, following Cartan's approach to generalized
geometry there is a minimal (i.e.\ locally described by a one dimensional
family of symplectic leaves) model
in this particular case: the contact projective geometry,
characterized by the maximally non-integrable symplectic distribution
of codimension one. This program was initiated in \cite{Kadlcakova2001},
and it immediately raised the question of a classification of
flat (homogeneous) invariants for (irreducible metaplectic
submodules of) the Segal--Shale--Weil
representation twisted by a one-parameter family of characters.

Relying on the recently developed approach to the
classification of homomorphisms between generalized
Verma modules called F-method,
\cite{koss}, \cite{Kostant1975}, \cite{Krizka-Somberg2017},
the present article shows a convenient setting given by
extending the category of induced modules to smooth admissible
representations.
As an example we give a complete
classification result in the case of
irreducible metaplectic submodules of the Segal--Shale--Weil
representation twisted by a family of characters as the inducing
representations, thereby supplying the missing structural information
(well known for finite-dimensional inducing representations.)
We recall that the idea of the F-method is to characterize the vector
space of all singular vectors in a generalized Verma module by a system of partial
differential equations acting on its geometric realization. The detailed explanation
of this idea is given in Section \ref{sec:Verma modules}. We also remark that for
characters as inducing representations
in the case of the $C_n$-series of simple finite-dimensional
Lie algebras, the composition structure for scalar generalized Verma
modules is rather poor (as we shall discuss in the forthcoming
article \cite{Krizka-Somberg2016})
when compared to the classical series $A_n,B_n,D_n$, see \cite{Krizka-Somberg2017}.
The vector-valued representation given by the
irreducible metaplectic submodules of the
Segal--Shale--Weil representation is a prominent example from the perspective of
potential applications.

We emphasize that the restriction to smooth inducing modules stems
from our main reason for the applicability, namely the construction
of equivariant differential operators on generalized flag manifolds.
Without this purpose, the
F-method itself can be treated for any inducing representation.

The content of our article goes as follows. In Section \ref{sec:Verma modules}, we briefly
introduce the formalism of the F-method (cf.\ \cite{koss}, \cite{Krizka-Somberg2017})
and extend it to a convenient category of smooth admissible inducing modules covering the finite-dimensional
representations. For this purpose we use a few standard arguments related to smooth globalization of
infinite-dimensional representations, cf.\ \cite{Bernstein-Krotz2014}, which in the same way
as for finite-dimensional representations
allows us to reformulate the task of finding geometrical invariants in the analytical framework
of solving a system of PDEs.
In Section \ref{sec:projective structure}, we describe the case of contact projective geometry in detail,
and introduce the family of representations induced from the irreducible metaplectic
submodules of the Segal--Shale--Weil representation twisted by characters. It is described in the non-compact
model on the open Schubert cell of the generalized flag manifold given by
the projective space. Section \ref{sec:singular vectors} presents
in Theorem \ref{thm:singular vector n}
a complete list of solutions of the former
problem in terms of singular vectors for the inducing
representations of our interest.
In Section \ref{sec:invariant operators general. flag}, we prove the
duality between equivariant differential
operators acting on induced representations and homomorphisms of
generalized Verma modules for smooth irreducible inducing representations,
a result well-known for finite-dimensional inducing
representations, cf.\ \cite{Collingwood-Shelton1990},
\cite{Soergel1990}.
Then we dualize the (inducing) representations and in Theorem
\ref{thm:edo} construct, following the complete set of singular vectors, all
equivariant differential operators on the irreducible metaplectic
submodules of the Segal--Shale--Weil representation
twisted by a character. Roughly speaking, our result
produces two types of equivariant differential operators. The first class
are the twistor-like operators landing in more complicated
metaplectic modules, and we expect their appearance in a not-yet-developed
Bernstein--Gelfand--Gelfand sequences associated
with infinite-dimensional representations. The second class
of equivariant operators is valued in the irreducible metaplectic
submodules of the Segal--Shale--Weil representation twisted by a character:
the operators have a right to be called
``contact powers of the contact symplectic Dirac operator'', in complete
analogy with the equivariant differential operators termed ``conformal
powers of the Laplace operator'', ``conformal
powers of the Dirac operator'' and ``CR invariant powers of the sub-Laplacian''
studied in the realm of the conformal and CR geometries, respectively,
cf.\ \cite{Graham-Jenne-Mason-Sparling1992}, \cite{Holland-Sparling2001}, \cite{Gover-Graham2005}.
The existence of first order invariant differential operators for the contact projective geometry and
the Segal--Shale--Weil representation was proved in \cite{Krysl2008}.

Throughout the article, the sets $\N$ and $\N_0$ denote $\{1,2,\dots\}$ and $\{0,1,2,\dots\}$, respectively.


\section{The structure of generalized Verma modules and equivariant differential operators}
\label{sec:Verma modules}

We consider the pair $(G,P)$ consisting of a connected real semisimple Lie group $G$ and
its parabolic subgroup $P$. In the Levi decomposition $P=LU$,
$L$ denotes the Levi subgroup and $U$ the unipotent radical of $P$.
We write $\gog(\mR),\,\gop(\mR),\,\gol(\mR),\,\mfrak{u}(\mR)$ for the
real Lie algebras and $\gog,\,\gop,\,\gol,\,\mfrak{u}$ for the complexified Lie
algebras of $G,\,P,\,L,\,U$, respectively. The symbols $U$ and $S$ applied to a Lie
algebra denote its universal enveloping algebra and symmetric algebra,
respectively.

It is well-known that the $G$-equivariant differential operators acting on
principal series representations for $G$ can be recognized in the study
of homomorphisms between generalized Verma modules for the Lie algebra
$\mathfrak{g}$. In the classical case of finite-dimensional inducing
representations, the latter homomorphisms are determined by
the image of inducing representation and its vectors are referred to as the
singular vectors. They are characterized as the vectors in generalized Verma
module annihilated by the positive nilradical $\mfrak{u}$.

As a starting point we use the F-method in order to classify and find precise positions of singular
vectors in a given representation space, see \cite{Kobayashi-Pevzner2016a}, \cite{koss}, \cite{Krizka-Somberg2017}
for a detailed exposition. In the present paper we generalize this technique
in the sense that we prove its natural extension from the category of finite-dimensional
inducing representations to a wider class of smooth admissible
inducing representations.
We briefly recall that a representation $(\sigma,\mV)$ of $P$ on a Hausdorff complete
locally convex topological (real or complex) vector space $\mV$
is smooth if all vectors in $\mV$ are smooth, and $v\in\mV$ is a smooth vector provided
the mapping from $P$ to $\mV$ given by $p \mapsto \sigma(p)v$ is smooth. A $U(\mfrak{p})$-module $\mV$ is admissible
if it is $\mfrak{Z}(\mfrak{l})$-finite, where $\mfrak{Z}(\mfrak{l})$ is the center of $U(\mfrak{l})$
and $\mfrak{Z}(\mfrak{l})$-finite means that the annihilator ideal of $\mV$ in
$\mfrak{Z}(\mfrak{l})$ is of finite codimension.
We refer to Section \ref{sec:invariant operators general. flag} for a more detailed discussion of
these notions.
\medskip

Let $(\sigma,\mathbb{V})$ be a smooth admissible $P$-module.
Let us assume that $\lambda \in \Hom_P(\mfrak{p},\C)$ defines a group character
$e^\lambda \colon P \rarr \GL(1,\C)$ of $P$ and define $\rho \in \Hom_P(\mfrak{p},\C)$ by
\begin{align}
  \rho(X) = {\textstyle {1 \over 2}} \tr_\mfrak{u}\ad(X) \label{eq:rho vector}
\end{align}
for all $X \in \mfrak{p}$, where $\tr_\mfrak{u} \ad(X)$ denotes the trace of $\ad(X) \colon \mfrak{u} \rarr \mfrak{u}$. Then we can define a twisted $P$-module $(\sigma_{\lambda+\rho},\mathbb{V}_{\lambda+\rho})$ with a twist $\lambda +\rho \in \Hom_P(\mfrak{p},\C)$ such that $p \in P$ acts as $e^{\lambda+\rho}(p)\sigma(p)v$ instead of
$\sigma(p)v$ for all $v \in \mathbb{V}_{\lambda+\rho} \simeq \mathbb{V}$ (an isomorphism of vector spaces).
In the case discussed by our main example, $\mV=\mathbb{S}$ is the complex smooth admissible simple
(unitarizable highest weight) $L$-module, $L \simeq \GL(1,\R) \times \Mp(2n,\R)$
with $n \in \N$, given by one of the metaplectic submodules in the
Segal--Shale--Weil representation.

For a chosen smooth principal series representation of $G$ on the vector space
$\Ind_P^G(\mV_{\lambda+\rho})$ of smooth sections of the homogeneous vector bundle
$G \times_P \mV_{\lambda+\rho} \rarr G/P$ associated to a complex smooth admissible
$P$-module
$\mV_{\lambda+\rho}$ (cf.\ \cite{Kriegl-Michor1997}), due to the smoothness of $\mV$
we obtain the infinitesimal action
\begin{align}
\pi_\lambda \colon \mfrak{g} \rarr \mcal{D}(U_e) \otimes_\C \End \mathbb{V}_{\lambda+\rho}
\end{align}
on the vector space $\mcal{C}^\infty(U_e) \otimes_\C\! \mathbb{V}_{\lambda+\rho}$
of $\mathbb{V}_{\lambda+\rho}$-valued smooth functions on $U_e$ in the non-compact picture of the induced representation.
Here $\mcal{D}(U_e)$ denotes the $\C$-algebra of smooth complex linear differential operators on $U_e=\widebar{U}P \subset G/P$ ($\widebar{U}$ is the Lie group whose Lie algebra is the opposite nilradical $\widebar{\mathfrak{u}}(\R)$ to $\mathfrak{u}(\R)$).

Since the vector space $\mcal{D}'_o(U_e)\otimes_\C\! \mV_{\lambda+\rho}$ of $\mathbb{V}_{\lambda+\rho}$-valued distributions on $U_e$ supported in the unit coset
$o=eP\in G/P$ is a $\mcal{D}(U_e) \otimes_\C \End \mathbb{V}_{\lambda+\rho}$-module, we get also
the infinitesimal action of $\pi_\lambda(X)$ for $X \in \mfrak{g}$ on $\mcal{D}'_o(U_e)\otimes_\C\! \mV_{\lambda+\rho}$. The exponential map allows us to identify $U_e$ with the nilpotent Lie algebra
$\widebar{\mfrak{u}}(\R)$. If we denote by $\eus{A}^\mfrak{g}_{\widebar{\mfrak{u}}}$
the Weyl algebra of the complex vector space $\widebar{\mfrak{u}}$, then the vector
space $\mcal{D}'_o(U_e)$ can be conveniently analyzed by identifying it as an
$\eus{A}^\mfrak{g}_{\widebar{\mfrak{u}}}$-module with the quotient of
$\eus{A}^\mfrak{g}_{\widebar{\mfrak{u}}}$ by the left ideal $I_e$ generated
by all polynomials on $\widebar{\mathfrak{u}}$ vanishing at the origin.
Moreover, there is a $U(\gog)$-module isomorphism
\begin{align}\label{vermadistr}
\Phi_\lambda \colon M^\gog_\gop(\mV_{\lambda-\rho})\equiv U(\gog)\otimes_{U(\gop)}\!\mV_{\lambda-\rho} \riso  \mcal{D}'_o(U_e) \otimes_\C\! \mV_{\lambda+\rho} \riso \eus{A}^\mfrak{g}_{\widebar{\mfrak{u}}}/I_e \otimes_\C\! \mathbb{V}_{\lambda+\rho}.
\end{align}
Let $(x_1,x_2,\dots,x_n)$ be the linear coordinate functions on $\widebar{\mfrak{u}}$ and let $(y_1,y_2,\dots,y_n)$ be the dual linear coordinate functions on $\widebar{\mfrak{u}}^*$. Then the algebraic Fourier transform
\begin{align}
  \mcal{F} \colon \eus{A}^\mfrak{g}_{\widebar{\mfrak{u}}} \rarr \eus{A}^\mfrak{g}_{\widebar{\mfrak{u}}^*}
\end{align}
is given by
\begin{align}
  \mcal{F}(x_i) = -\partial_{y_i}, \qquad \mcal{F}(\partial_{x_i}) = y_i
\end{align}
for $i=1,2,\dots,n$, and leads to a vector space isomorphism
\begin{align}\label{eqn:FT}
\tau \colon  \eus{A}^\mfrak{g}_{\widebar{\mfrak{u}}}/I_e  \riso  \eus{A}^\mfrak{g}_{\widebar{\mfrak{u}}^*}/\mcal{F}(I_e) \simeq \C[\widebar{\mfrak{u}}^*]
\end{align}
defined by
\begin{align}
  Q\ {\rm mod}\ I_e \mapsto \mcal{F}(Q)\ {\rm mod}\ \mcal{F}(I_e)
\end{align}
for $Q \in \eus{A}^\mfrak{g}_{\widebar{\mfrak{u}}}$. The composition of the previous mappings \eqref{vermadistr} and \eqref{eqn:FT} gives the vector space isomorphism
\begin{align}\label{eqn:phiPM}
 \tau \circ \Phi_\lambda \colon U(\mfrak{g}) \otimes_{U(\mfrak{p})} \! \mathbb{V}_{\lambda-\rho} \riso \eus{A}^\mfrak{g}_{\widebar{\mfrak{u}}}/I_e \otimes_\C \!\mV_{\lambda-\rho} \riso
\C[\widebar{\mathfrak{u}}^*]\otimes_\C\! \mV_{\lambda+\rho} \riso \C[\widebar{\mathfrak{u}}^*]\otimes_\C\! \mV_{\lambda-\rho},
\end{align}
where the last conventional isomorphism accounts for the identification between the left and right
$\mcal{D}$-module structure, and thereby it induces an action
\begin{align}
\hat{\pi}_\lambda \colon \mfrak{g} \rarr \eus{A}^\mfrak{g}_{\widebar{\mfrak{u}}^*}\! \otimes_\C \End \mathbb{V}_{\lambda-\rho}
\end{align}
of $\mfrak{g}$ on $\C[\widebar{\mathfrak{u}}^*]\otimes_\C\! \mV_{\lambda-\rho}$.

In order to find an explicit form of $G$-equivariant differential operators, discussed in more detail in Section \ref{sec:invariant operators general. flag}, it is necessary to have an explicit formula for the inverse of the mapping $\tau \circ \Phi_{\lambda-\rho} \colon M^\mfrak{g}_\mfrak{p}(\mathbb{V}_{\lambda-\rho}) \rarr \C[\widebar{\mathfrak{u}}^*]\otimes_\C\! \mathbb{V}_{\lambda-\rho}$. Let us introduce the symmetrization map $\beta \colon S(\widebar{\mfrak{u}}) \rarr U(\widebar{\mfrak{u}})$, which is a $\C$-linear isomorphism of filtered vector spaces defined on monomials by
\begin{align}
  \beta(a_1a_2\dots a_k) = {1 \over k!} \sum_{\sigma \in S_k} a_{\sigma(1)}
	a_{\sigma(2)} \dots a_{\sigma(k)}  \label{symmetr}
\end{align}
for all $k\in \N$ and $a_1,a_2,\dots,a_k \in \widebar{\mfrak{u}}$, where $S_k$ is the permutation group of the set $\{1,2,\dots,k\}$.
If we denote by $(f_1,f_2,\dots,f_n)$ a basis of $\widebar{\mfrak{u}}$ satisfying $\hat{\pi}_\lambda(f_i)=-y_i\ {\rm mod}\ \mcal{F}(I_e)$ for $i=1,2,\dots,n$, then we have
\begin{align}
(\tau \circ \Phi_\lambda)(\beta(f_{i_1}f_{i_2}\dots f_{i_k}) \otimes v)= (-1)^k y_{i_1}y_{i_2}\dots y_{i_k} \otimes v  \label{eq:inverse symmetrization}
\end{align}
for all $v \in \mathbb{V}_{\lambda-\rho}$, $k \in \N$ and $i_1,i_2,\dots,i_k \in \{1,2,\dots,n\}$.

\medskip

\definition{\label{vmgendef}Let $\mV$ be a smooth admissible $L$-module which extends to a
$P$-module by $U$ acting trivially. Then the generalized Verma module
induced from $\mV$ is the $(\mfrak{g},P)$-module
\begin{align}\label{vmdefinition}
  M^\mfrak{g}_\mfrak{p}(\mathbb{V}) = \Ind^\mfrak{g}_\mfrak{p} \mathbb{V} \equiv U(\mfrak{g}) \otimes_{U(\mfrak{p})}\!\mathbb{V}  \simeq U(\widebar{\mfrak{u}}) \otimes_\C\!\mathbb{V},
\end{align}
where the last isomorphism of $U(\widebar{\mfrak{u}})$-modules follows from the Poincaré--Birkhoff--Witt theorem.}

By abuse of notation we call the modules $M^\mfrak{g}_\mfrak{p}(\mathbb{V})$ generalized Verma modules,
although this terminology is usually reserved for modules induced from finite-dimensional representations
(cf.\ \cite{Mazorchuk-Stroppel2008} for the terminology in the case of inducing ${\mfrak l}$-modules
in the BGG category ${\fam2 O}$ associated to $U({\mfrak l})$). The generalized Verma module
$M^\mfrak{g}_\mfrak{p}(\mathbb{V})$ has a $P$-module structure defined as follows. As
$M^\mfrak{g}_\mfrak{p}(\mathbb{V})$ is locally $\mfrak{u}$-finite and $U$ is simply connected,
a $U$-module structure on $M^\mfrak{g}_\mfrak{p}(\mathbb{V})$ is given by the exponential
mapping $\exp \colon \mfrak{u} \rarr U$. Since $U(\widebar{\mfrak{u}})$ is an $L$-module by
the adjoint action of $L$ and $\mathbb{V}$ is by definition a $L$-module, we have a $L$-module
structure on $U(\widebar{\mfrak{u}}) \otimes_\C \!\mathbb{V}$.
The Levi decomposition $P=LU$ then results in a $P$-module structure on $M^\mfrak{g}_\mfrak{p}(\mathbb{V})$.

As $U(\widebar{\mfrak{u}})$ is a direct sum of finite-dimensional $L$-modules, the tensor product
$U(\widebar{\mfrak{u}}) \otimes_\C \!\mathbb{V}$ is a direct sum of the tensor products of
finite-dimensional $L$-modules with ${\mathbb{V}}$. We notice that basic properties of $L$-module composition
structure of such tensor products can be found for example in \cite{Bernstein-Gelfand1980}.

Based on the language of globalizations of infinite-dimensional representations it is well-known that
the tensor product of a smooth module (smooth globalization of a Harish-Chandra module) with a
finite-dimensional representation is a smooth module, and a direct sum of smooth modules is a smooth module.
Moreover, quotients of smooth modules by smooth submodules are smooth modules, cf.\ \cite{Bernstein-Krotz2014}
for all these results. In particular, $M^\mfrak{g}_\mfrak{p}(\mathbb{V})$ is a smooth $(\mfrak{g},P)$-module
which is a direct sum of smooth admissible $L$-modules because $U(\widebar{\mfrak{u}})$ is a direct sum of
finite-dimensional $L$-modules.

\medskip

\definition{Let $\mV$ be a smooth admissible $L$-module which extends to a
$P$-module by $U$ acting trivially. We define an $L$-module
\begin{align}
M_\mfrak{p}^\mfrak{g}(\mV)^\mathfrak{u}= \{v\in M^\mfrak{g}_\mfrak{p}(\mV);\, Xv=0\
\text{for all}\  X\in \mathfrak{u}\},
\end{align}
and call it the vector space of singular vectors.}

Let us note that the $L$-module structure on $M^\mfrak{g}_\mfrak{p}(\mathbb{V})^\mfrak{u}$ is induced by the $P$-module structure on the generalized Verma module $M^\mfrak{g}_\mfrak{p}(\mathbb{V})$. Since $\Ad(l^{-1})(X) \in \mfrak{u}$ for all $X \in \mfrak{u}$ and $l \in L$, where $\Ad(l^{-1})$ denotes the adjoint action on $\mfrak{g}$, we obtain $X(lv)=l\Ad(l^{-1})(X)v=0$ for all $X \in \mfrak{u}$, $l \in L$, $v \in M^\mfrak{g}_\mfrak{p}(\mathbb{V})^\mfrak{u}$ and so
$lv \in M^\mfrak{g}_\mfrak{p}(\mathbb{V})^\mfrak{u}$ for all $l \in L$, $v \in M^\mfrak{g}_\mfrak{p}(\mathbb{V})^\mfrak{u}$.

We also notice that for a pair of smooth admissible $L$-modules $\mathbb{W}$ and $\mathbb{V}$, any
$(\mfrak{g},P)$-module homomorphism
\begin{align}
\varphi \colon  M^\mfrak{g}_\mfrak{p}(\mathbb{W})\to M^\mfrak{g}_\mfrak{p}(\mathbb{V})
\end{align}
is determined by the image $\varphi(1\otimes \mathbb{W})$, because
$M^\mfrak{g}_\mfrak{p}(\mathbb{W})$ is a free $U(\widebar{\mfrak{u}})$-module by \eqref{vmdefinition}.
\medskip

The main difference between \cite{koss} and its generalization considered here is that
$M_\mfrak{p}^\mfrak{g}(\mV)^\mathfrak{u}$ does not decompose into direct sum of simple
$L$-modules in general.
Let us denote by $\mW$ an $L$-submodule of $M_\mfrak{p}^\mfrak{g}(\mV)^\mathfrak{u}$. By the characterizing
property of $(\mfrak{g},P)$-module homomorphism in the category of modules
$M^\mfrak{g}_\mfrak{p}(\mathbb{V})$ mentioned in the previous paragraph, we have
\begin{align}
\Hom_{(\mfrak{g},P)}(M_\mfrak{p}^\mfrak{g}(\mW), M_\mfrak{p}^\mfrak{g}(\mV))
\simeq
\Hom_L(\mW, M_\mfrak{p}^\mfrak{g}(\mV)^\mathfrak{u}),
\end{align}
where $M_\mfrak{p}^\mfrak{g}(\mW)$ has the composition series given by smooth admissible $(\mfrak{g},P)$-modules.

In the main application treated in Section \ref{sec:projective structure} and Section \ref{sec:singular vectors},
we apply the previous considerations to the infinite-dimensional smooth (highest weight unitarizable) $L$-module
given by the simple $L$-submodules of the Segal--Shale--Weil representation, see \cite{DeBie-Somberg-Soucek2014}.
The space of singular vectors has in this case an unexpected feature, namely, it is a direct sum of
infinite-dimensional simple highest weight $L$-modules.

The rest of the F-method as present in \cite{koss} works without any modification for smooth admissible
inducing modules $\mathbb{V}$. In particular, the isomorphism
$\tau \circ \Phi_{\lambda+\rho} \colon M^\mfrak{g}_\mfrak{p}(\mathbb{V}_{\lambda-\rho}) \rarr \C[\widebar{\mfrak{u}}^*] \otimes_\C\! \mathbb{V}_{\lambda-\rho}$ of $\mfrak{g}$-modules ensures that $\C[\widebar{\mfrak{u}}^*] \otimes_\C\! \mathbb{V}_{\lambda-\rho}$ is a $(\mfrak{g},P)$-module as well. We introduce an $L$-module
\begin{align} \label{eqn:sol2l}
\Sol(\mathfrak{g},\mathfrak{p};\C[\widebar{\mathfrak{u}}^*]\otimes_\C\! \mV_{\lambda-\rho})^\mcal{F} =\{f \in \C[\widebar{\mathfrak{u}}^*]\otimes_\C\! \mV_{\lambda-\rho};\, \hat{\pi}_\lambda(X) f = 0\ \text{for all}\ X\in\mathfrak{u}\},
\end{align}
which is in fact the vector space of singular vectors in $\C[\widebar{\mfrak{u}}^*] \otimes_\C\! \mathbb{V}_{\lambda-\rho}$,
and by \eqref{eqn:phiPM} there is an $L$-equivariant isomorphism
\begin{equation}
\label{eqn:phi}
\tau \circ \Phi_\lambda \colon M_\mathfrak{p}^\mathfrak{g}(\mV_{\lambda-\rho})^\mathfrak{u} \riso \Sol(\mathfrak{g},\mathfrak{p}; \C[\widebar{\mathfrak{u}}^*]\otimes_\C\! \mV_{\lambda-\rho})^\mcal{F}.
\end{equation}
The action of $\hat{\pi}_\lambda(X)$ on $\C[\widebar{\mathfrak{u}}^*]\otimes_\C\! \mV_{\lambda-\rho}$
produces a system of partial differential equations for the elements in
$\Sol(\mathfrak{g},\mathfrak{p}; \C[\widebar{\mathfrak{u}}^*]\otimes_\C\! \mV_{\lambda-\rho})^\mcal{F}$, which makes it possible to describe its structure completely in particular cases of interest. Namely, the algebraic Fourier transform on Weyl algebra modules converts the algebraic problem of finding singular vectors in generalized Verma
modules into an analytic problem of solving the systems of partial differential equations.

The formulation above has a classical dual statement
(cf.\ \cite{Collingwood-Shelton1990}, \cite{Soergel1990} for the
formulation in the category of finite-dimensional inducing $P$-modules), which
explains the relationship between the geometrical problem of finding
$G$-equivariant differential operators between induced representations
and the algebraic problem of finding homomorphisms between generalized
Verma modules. We refer to Section \ref{sec:invariant operators general. flag}
for a detailed exposition of this duality for an arbitrary inducing
smooth admissible $P$-module.

\section{The contact projective structure}
\label{sec:projective structure}

In the present section we describe a few basic geometrical
and representation theoretical aspects of the flat (or homogeneous)
contact projective structure needed in our analysis. For the purposes
of applications, we review the real homogeneous model of
the contact projective geometry.

\subsection{The geometry of flat contact projective structure}

The generalized flag manifold corresponding to the flat model of
real contact projective structure is the homogeneous space
$G/P\simeq \mathbb{RP}^{2n+1}$, where
the group of symplectic automorphisms $G=\Sp(2n+2,\mR)$ of the standard
symplectic vector space $(\mR^{2n+2},\Omega)$, $n\in\mN$, acts transitively
on the space of lines in
$\mR^{2n+2}$ by $(g,[v])\mapsto [g.v]$ for $0\not=v\in\mR^{2n+2}$,
$g\in \Sp(2n+2,\mR)$. The stabilizer of a line is a parabolic subgroup
conjugate to the standard parabolic subgroup
$P=(\GL(1,\R) \times \Sp(2n,\mR))\ltimes {\rm H}(n,\R)$, where ${\rm H}(n,\R)$ is the
real Heisenberg group of dimension $2n+1$. In the
Dynkin diagrammatic notation, the Lie algebra of $P$ is given
by omitting the first simple root in the $C_{n+1}$-series of simple Lie
algebras.

The geometry of real contact projective structure can be described as follows.
We denote by $V^1\subset \mR^{2n+2}$ the $1$-dimensional subspace whose
stabilizer is the group $P$, so that its $\Omega$-complement $V^2=(V^1)^\perp$
yields the $P$-invariant filtration $V^1\subset V^2\subset\mR^{2n+2}$. The choice
$0\not= v_\infty\in V^1$ and $v_0\in\mR^{2n+2}$ such that
$\Omega(v_\infty,v_0)=1$ determines a splitting
$\mR^{2n+2}= V_0\oplus V_{-1}\oplus V_{-2}$ with $V_0=V^1$
and $v_0$ spans $V_{-2}$.
This grading on $\mR^{2n+2}$ induces a
$\mZ$-grading on the Lie algebra $\gog(\mR)$ of $G$,
$\gog(\mR)=\oplus_{i=-2}^2\gog(\mR)_i$, by
\begin{align*}
\gog(\mR)_i=\{X\in\gog(\mR);\, X(V_k)\subset V_{k+i}\ \text{for}\ -2\leq k\leq 0 \}.
\end{align*}
The subalgebra $\gog(\mR)_{-2}\oplus\gog(\mR)_{-1}$ is the Heisenberg algebra,
$\dim \gog(\mR)_{-2}=1$, $\dim \gog(\mR)_{-1}=2n$, and $\gog(\mR)_0$ is the conformal
symplectic Lie algebra. The vectors $v_0$ and $v_\infty$ together with a basis
$v_1,\dots ,v_{2n}$ of $V_{-1}$ then determine linear coordinate functions $x^0,x^\infty,x^1,\dots, x^{2n}$
on $\mR^{2n+2}$, such that $\Omega=\frac{1}{2}\Omega_{IJ}dx^I\wedge dx^J$,
$I,J\in\{0,\infty,1,\dots ,2n\}$, with
\begin{align*}
\Omega_{0\infty}=1,\qquad \Omega_{0i}=0,\qquad \Omega_{\infty i}=0, \qquad \Omega_{ij}=\omega_{ij},
\end{align*}
where $\omega_{ij}=\Omega(v_i,v_j)$. The group of automorphisms of $(\mR^{2n+2},\Omega)$ is $G$, the Lie algebra
of $P$ is $\gog(\mR)_0\oplus\gog(\mR)_1\oplus\gog(\mR)_2$. We denote the stabilizer of the vector
$v_\infty$ by $P^s\subset P$, its Lie algebra by $\gop(\mR)^s$.
Since we have $P/P^s \simeq \GL(1,\R)$, the principal $\GL(1,\R)$-bundle $G/P^s \rarr G/P$ determines the tautological vector bundle
$\mR^{2n+2}\setminus\{0\}\rarr  \mR\mP^{2n+1}$, and the Maurer-Cartan form on $G$
is the Cartan connection for the flat contact projective structure.

The $P$-invariant filtration on $\gog(\mR)$ allows us to define on $\gog(\mR)_{-1}$
the conformal symplectic structure, which further determines a symplectic structure on
$\gog(\mR)/\gop(\mR)^s$ isomorphic to $(\mR^{2n+2},\Omega)$. The basis
$v_0,v_1,\dots ,v_{2n}$ of $\gog(\mR)_{-2}\oplus\gog(\mR)_{-1}$ can be exponentiated
to the basis of left-invariant vector fields $E_0,E_1,\dots ,E_{2n}$ on $G/P$, with
the Lie bracket $[E_\alpha,E_\beta]=-\Omega_{\alpha\beta} E_0$,
$\alpha,\beta \in\{0,1,\dots ,2n\}$. The canonical
contact structure (or, the contact distribution) on $T(G/P)$ is the
left $G$-invariant subbundle of $T(G/P)$,
generated by $\gog(\mR)_{-1}$ and spanned at each point by $E_1,\dots ,E_{2n}$.
A left-invariant section $\theta\in T^*(G/P)$ of the annihilator of the contact
distribution is
then determined by the requirement that $E_0$ is its Reeb vector field. The
covariant derivative $\nabla$, defined by the requirement that the vector
fields $E_\alpha$ are $\nabla$-parallel, preserves $\theta$ and $d\theta$ and
determines a model for any contact projective geometry.

From the representation-theoretical point of view, we need the double cover
of the symplectic group $G=\Sp(2n+2,\R)$ called the metaplectic group
$\smash{\widetilde{G}}=\Mp(2n+2,\R)$. Thus we have a double cover homomorphism
$\psi  \colon \smash{\widetilde{G}} \rarr G$ of Lie groups, and
if we define the parabolic subgroup $\smash{\widetilde{P}}$ of
$\smash{\widetilde{G}}$ by $\smash{\widetilde{P}}=\psi^{-1}(P)$,
then the generalized flag manifold $\smash{\widetilde{G}/\widetilde{P}}$
is isomorphic to $G/P$ induced by the mapping $\psi$.

For more detailed exposition with a view towards the description of a
general (curved) contact projective structure, cf.\ \cite{Fox2005} and the
references therein.

\subsection{Representation theory of contact projective structure}
\label{section3.1}

Let us consider the connected complex simple Lie group $G_\C=\Sp(2n+2,\C)$, $n \geq 1$, defined by
\begin{align}
  \Sp(2n+2,\C)=\{g\in \GL(2n+2,\C);\, g^{\rm T}\!J_{2n+2}g=J_{2n+2}\}, \quad
  J_{2n+2}=\begin{pmatrix}
    0 & I_{n+1} \\
    -I_{n+1} & 0
  \end{pmatrix}\!,
\end{align}
and its Lie algebra $\mfrak{g}=\mfrak{sp}(2n+2,\C)$ given by
\begin{align}
\begin{aligned}
 \mfrak{sp}(2n+2,\C)&=\{X\in M_{2n+2,2n+2}(\C);\, X^{\rm T}\!J_{2n+2}+J_{2n+2}X=0\}  \\
& =\left\{\!
  \begin{pmatrix}
    A & B \\
    C & -A^{\rm T}
  \end{pmatrix}\!;\,
  A,B,C\in M_{n+1,n+1}(\C),\, B^{\rm T}=B,\, C^{\rm T}=C\right\}\!.
\end{aligned}
\end{align}
A Cartan subalgebra $\mfrak{h}$ of $\mfrak{g}$ is defined by the diagonal matrices
\begin{align}
  \mfrak{h}=\{\diag(a_1, \dots, a_{n+1}, -a_1, \dots, -a_{n+1});\, a_1,a_2,\dots,a_{n+1} \in \C\}.
\end{align}
For $i=1,2,\dots, n+1$ we define $\veps_i \in \mfrak{h}^*$ by $\veps_i(\diag(a_1, \dots, a_{n+1}, -a_1, \dots, -a_{n+1}))=a_i$. Then the root system of $\mfrak{g}$ with respect to $\mfrak{h}$ is
$\Delta = \{\pm \veps_i \pm \veps_j;\, 1 \leq i < j \leq n+1\}\cup\{\pm 2\veps_i;\, 1 \leq i \leq n+1\}$. A positive root system is $\Delta^+=\{\veps_i\pm\veps_j;\, 1\leq i<j \leq n+1\}\cup\{2\veps_i;\, 1 \leq i \leq n+1\}$ with the set of simple roots $\Pi=\{\alpha_1,\alpha_2,\dots,\alpha_{n+1}\}$, $\alpha_i=\veps_i-\veps_{i+1}$, $i=1, 2, \dots, n$ and $\alpha_{n+1}=2\veps_{n+1}$. Finally, the fundamental weights are $\omega_i= \smash{{\textstyle \sum_{j=1}^i}} \veps_j, i=1,2,\dots,n+1$.

The Lie subalgebras $\mfrak{b}$ and $\widebar{\mfrak{b}}$ defined as the direct sum of positive and negative root spaces together with the Cartan subalgebra are called the standard Borel and the opposite standard Borel subalgebras of $\mfrak{g}$, respectively. The subset
$\Sigma=\{\alpha_2,\alpha_3,\dots,\alpha_{n+1}\}$ of $\Pi$ generates the root
subsystem $\Delta_\Sigma$ in $\mfrak{h}^*$, and we associate to $\Sigma$ the
standard parabolic subalgebra $\mfrak{p}$ of $\mfrak{g}$ by $\mfrak{p} = \mfrak{l} \oplus \mfrak{u}$.
The reductive Levi subalgebra $\mfrak{l}$ of $\mfrak{p}$ is defined through
\begin{align}
\mfrak{l}= \mfrak{h} \oplus \bigoplus_{\alpha \in \Delta_\Sigma} \mfrak{g}_\alpha,
\end{align}
and the nilradical $\mfrak{u}$ of $\mfrak{p}$ and the opposite nilradical $\widebar{\mfrak{u}}$ are
\begin{align}
  \mfrak{u}= \bigoplus_{\alpha \in \Delta^+ \setminus \Delta_\Sigma^+}\mfrak{g}_\alpha \qquad \text{and}
  \qquad \widebar{\mfrak{u}}= \bigoplus_{\alpha \in \Delta^+ \setminus \Delta_\Sigma^+} \mfrak{g}_{-\alpha},
\end{align}
respectively. We define the $\Sigma$-height $\htt_\Sigma(\alpha)$ of $\alpha \in \Delta$ by
\begin{align}
  \htt_\Sigma\!\big({\textstyle \sum_{i=1}^{n+1}} a_i \alpha_i\big) = a_1,
\end{align}
so $\mfrak{g}$ is a $|2|$-graded Lie algebra with respect to the grading given by $\mfrak{g}_i = \bigoplus_{\alpha \in \Delta,\, \htt_\Sigma(\alpha)=i} \mfrak{g}_\alpha$ for $0 \neq i \in \Z$, and $\mfrak{g}_0= \mfrak{h} \oplus \smash{\bigoplus}_{\alpha \in \Delta,\, \htt_\Sigma(\alpha)=0} \mfrak{g}_\alpha$. Moreover, we have $\mfrak{u}=\mfrak{g}_1 \oplus \mfrak{g}_2$, $\widebar{\mfrak{u}}=\mfrak{g}_{-2} \oplus \mfrak{g}_{-1}$ and $\mfrak{l}= \mfrak{g}_0 \simeq \C \oplus \mfrak{sp}(2n,\C)$.

Let us denote by $1_i$, $i=1,2,\dots,n$, the $(n \times 1)$-matrix having $1$ at the $i$-th row and $0$ elsewhere. The basis $(f_1,\dots,f_n,g_1,\dots,g_n,c)$ of the root spaces in the opposite nilradical $\widebar{\mfrak{u}}$ is given by
\begin{align}
  f_i = \begin{pmatrix}
    0 & 0 & 0 & 0 \\
    1_i & 0 & 0 & 0 \\
    0 & 0 & 0 & -1_i^{\rm T} \\
    0 & 0 & 0 & 0
  \end{pmatrix}\!, \qquad
  g_i = \begin{pmatrix}
    0 & 0 & 0 & 0 \\
    0 & 0 & 0 & 0 \\
    0 & 1_i^{\rm T} & 0 & 0 \\
    1_i & 0 & 0 & 0
  \end{pmatrix}\!, \qquad
  c = \begin{pmatrix}
    0 & 0 & 0 & 0 \\
    0 & 0 & 0 & 0 \\
    2 & 0 & 0 & 0 \\
    0 & 0 & 0 & 0
  \end{pmatrix}\!,
\end{align}
where the only non-trivial Lie brackets are $[f_i,g_i]=-c$ for all $i=1,2,\dots,n$. Analogously, the basis $(d_1,\dots,d_n,e_1,\dots,e_n,a)$ of the root spaces in $\mfrak{u}$ is given by
\begin{align}
  d_i = \begin{pmatrix}
    0 & 1_i^{\rm T} & 0 & 0 \\
    0 & 0 & 0 & 0 \\
    0 & 0 & 0 &  0 \\
    0 & 0 & -1_i & 0
  \end{pmatrix}\!, \qquad
  e_i = \begin{pmatrix}
    0 & 0 & 0 & 1_i^{\rm T} \\
    0 & 0 & 1_i & 0 \\
    0 & 0 & 0 & 0 \\
    0 & 0 & 0 & 0
  \end{pmatrix}\!, \qquad
  a = \begin{pmatrix}
    0 & 0 & 2 & 0 \\
    0 & 0 & 0 & 0 \\
    0 & 0 & 0 & 0 \\
    0 & 0 & 0 & 0
  \end{pmatrix}\!,
\end{align}
where $[d_i,e_i]=a$ for all $i=1,2,\dots,n$. The Levi subalgebra $\mfrak{l}$ of $\mfrak{p}$ is the linear span of
\begin{align}
  h = \begin{pmatrix}
    1 & 0 & 0 & 0 \\
    0 & 0 & 0 & 0 \\
    0 & 0 & -1 &  0 \\
    0 & 0 & 0 & 0
  \end{pmatrix}\!, \qquad
  h_{A,B,C} = \begin{pmatrix}
    0 & 0 & 0 & 0 \\
    0 & A & 0 & B \\
    0 & 0 & 0 &  0 \\
    0 & C & 0 & -A^{\rm T}
  \end{pmatrix}\!,
\end{align}
where $A,B,C \in M_{n,n}(\C)$ satisfy $B^{\rm T}=B$ and $C^{\rm T}=C$. Moreover,
the element $h$ is a basis of the center $\mfrak{z}(\mfrak{l})$ of $\mfrak{l}$.

Finally, the parabolic subgroup $P_\C$ of $G_\C$ with the Lie algebra $\mfrak{p}$ is defined by $P_\C=N_G(\mfrak{p})$.
The parabolic subalgebra $\mfrak{p} \subset \mfrak{g}$ is given by
\begin{align}
  \mfrak{p}=\left\{\begin{pmatrix}
    h & d^{\rm T} & a & e^{\rm T}\\
    0 & A & e & B \\
    0 & 0 & -h & 0 \\
    0 & C & -d & -A^{\rm T}
  \end{pmatrix}\!;
  \begin{gathered}
   h,a\in \C,\, d,e \in \C^n,\, A,B,C \in M_{n,n}(\C), \\
   B^{\rm T}=B,\, C^{\rm T}=C
   \end{gathered}
  \right\}\!.
\end{align}
The real connected simple Lie group $G$ and its parabolic subgroup $P$ is defined by $G=G_\C \cap \GL(2n+2,\R)$ and $P=P_\C \cap \GL(2n+2,\R)$.

Any character $\lambda \in \Hom_P(\mfrak{p},\C)$ is given by
\begin{align}
  \lambda= \lambda_1 \widetilde{\omega}_1
\end{align}
for some $\lambda_1\in \C$, where $\widetilde{\omega}_1\in \Hom_P(\mfrak{p},\C)$
is equal to $\omega_1\in \mfrak{h}^*$ regarded as trivially extended to
$\mfrak{p} =\mfrak{h} \oplus (\bigoplus_{\alpha \in \Delta_\Sigma} \mfrak{g}_\alpha) \oplus \mfrak{u}$.
The vector $\rho \in \Hom_P(\mfrak{p},\C)$ defined by \eqref{eq:rho vector} is
\begin{align}
  \rho=(n+1)\widetilde{\omega}_1.
\end{align}
By abuse of notation, for the character $\lambda\widetilde{\omega}_1 \in \Hom_P(\mfrak{p},\C)$, $\lambda \in \C$, we use the simplified notation $\lambda \in \Hom_P(\mfrak{p},\C)$.

\subsection{Description of the representation}

Here we describe the representations of $\mathfrak{g}$
on the space of sections of vector bundles on $G/P$ associated to the
simple metaplectic submodules of the Segal--Shale--Weil representation
$\mathbb{S}_{\lambda+\rho}$ of $\smash{\widetilde{P}}$ twisted by
characters $\lambda + \rho \in \Hom_P(\mfrak{p},\C)$.

The induced representations in question are described in
the non-compact picture, given by restricting sections on $G/P$ to the
open Schubert cell $U_e$ isomorphic by the exponential map to the opposite
nilradical $\widebar{\mfrak{u}}(\R)$.
Let us denote by $(\hat{x}_1,\dots,\hat{x}_n,\hat{y}_1,\dots,\hat{y}_n,\hat{z})$ the linear coordinate functions on
$\widebar{\mfrak{u}}$ with respect to the basis $(f_1,\dots,f_n,g_1,\dots,g_n,c)$ of the opposite nilradical $\widebar{\mfrak{u}}$, and by $(x_1,\dots,x_n,y_1,\dots,y_n,z)$ the dual linear coordinate functions on $\widebar{\mfrak{u}}^*$. Then the Weyl algebra
$\eus{A}^\mfrak{g}_{\widebar{\mfrak{u}}}$ is generated by
\begin{align}
\{\hat{x}_1,\dots,\hat{x}_n,\hat{y}_1,\dots,\hat{y}_n,\hat{z},\partial_{\hat{x}_1},\dots,\partial_{\hat{x}_n}, \partial_{\hat{y}_1},\dots, \partial_{\hat{y}_n},\partial_{\hat{z}}\}
\label{eq:Weyl algebra generators}
\end{align}
and the Weyl algebra
$\eus{A}^\mfrak{g}_{\widebar{\mfrak{u}}^*}$ is generated by
\begin{align}
\{x_1,\dots,x_n,y_1,\dots,y_n,z,\partial_{x_1},\dots,\partial_{x_n}, \partial_{y_1},\dots, \partial_{y_n},\partial_{z}\}.
\label{eq:dual Weyl algebra generators}
\end{align}
The local coordinate chart $u_e \colon x\in U_e \mapsto u_e(x)\in \widebar{\mfrak{u}}(\R) \subset \widebar{\mfrak{u}}$ for the open subset $U_e\subset G/P$, in coordinates with respect to the basis $(f_1,\dots,f_n,g_1,\dots,g_n,c)$ of $\widebar{\mfrak{u}}$, is given by
\begin{align}
  u_e(x)=\sum_{i=1}^n u^i_x(x)f_i + \sum_{i=1}^n u^i_y(x) g_i + u_z(x)c \label{eq:coordinate function}
\end{align}
for all $x \in U_e$.
\medskip

Let $(\sigma,\mathbb{V})$, $\sigma \colon \mfrak{p} \rarr \mfrak{gl}(\mathbb{V})$,
be a $\mfrak{p}$-module. Then a twisted $\mfrak{p}$-module $(\sigma_\lambda, \mathbb{V}_\lambda)$,
$\sigma_\lambda \colon \mfrak{p} \rarr \mfrak{gl}(\mathbb{V}_\lambda)$, with a twist
$\lambda \in \Hom_P(\mfrak{p},\C)$, is defined by
\begin{align}
  \sigma_\lambda(X)v=\sigma(X)v+\lambda(X)v
\end{align}
for all $X \in \mfrak{p}$ and $v \in \mathbb{V}_\lambda \simeq \mathbb{V}$ (as vector spaces).
\medskip

Let us introduce the notation
\begin{align}
  E_x = \sum_{j=1}^n x_j\partial_{x_j}, \qquad E_z=z\partial_z, \qquad E_y = \sum_{j=1}^n y_j\partial_{y_j}
\end{align}
and
\begin{align}
  E_{\hat{x}} = \sum_{j=1}^n \hat{x}_j\partial_{\hat{x}_j}, \qquad E_{\hat{z}}=\hat{z}\partial_{\hat{z}}, \qquad E_{\hat{y}} = \sum_{j=1}^n \hat{y}_j\partial_{\hat{y}_j}
\end{align}
for the Euler homogeneity operators.
\medskip

\theorem{\label{reproper}
Let $\lambda \in \Hom_P(\mfrak{p},\C)$ and let $(\sigma,\mathbb{V})$ be a $\mfrak{p}$-module. Then the embedding of $\mfrak{g}$ into
$\eus{A}^\mfrak{g}_{\widebar{\mfrak{u}}} \otimes_\C \End \mathbb{V}_{\lambda+\rho}$ and
$\eus{A}^\mfrak{g}_{\widebar{\mfrak{u}}^*}\! \otimes_\C \End \mathbb{V}_{\lambda-\rho}$ is given by
\begin{enumerate}
\item[1)]
\begin{align}
\begin{aligned}
  \pi_\lambda(f_i)&=-\partial_{\hat{x}_i}+{\textstyle {1\over 2}}\hat{y}_i\partial_{\hat{z}},  \\
  \pi_\lambda(g_i)&=-\partial_{\hat{y}_i}-{\textstyle {1\over 2}}\hat{x}_i\partial_{\hat{z}},  \\
  \pi_\lambda(c)&=-\partial_{\hat{z}},
\end{aligned}
\end{align}
\begin{align}
\begin{aligned}
  \hat{\pi}_\lambda(f_i)&=-x_i-{\textstyle {1\over 2}}z\partial_{y_i}, \\
  \hat{\pi}_\lambda(g_i)&=-y_i+{\textstyle {1\over 2}}z\partial_{x_i}, \\
  \hat{\pi}_\lambda(c)&=-z
\end{aligned}
\end{align}
for $i=1,2,\dots,n$;
\item[2)]
\begin{align}
\begin{aligned}
  \pi_\lambda(h)&= E_{\hat{x}}  + E_{\hat{y}}+ 2E_{\hat{z}} + \sigma_{\lambda+\rho}(h), \\
  \pi_\lambda(h_{A,B,C})&=- {\textstyle \sum\limits_{i,j=1}^n}
	a_{ij}(\hat{x}_j\partial_{\hat{x}_i}-\hat{y}_i\partial_{\hat{y}_j})
	-{\textstyle \sum\limits_{i,j=1}^n} (b_{ij}\hat{y}_j\partial_{\hat{x}_i}+c_{ij}\hat{x}_j\partial_{\hat{y}_i}) +\sigma_{\lambda+\rho}(h_{A,B,C}),
\end{aligned}
\end{align}
\begin{align}
\begin{aligned}
  \hat{\pi}_\lambda(h)&= -E_x  - E_y - 2E_z + \sigma_{\lambda-\rho}(h), \\
  \hat{\pi}_\lambda(h_{A,B,C})&= {\textstyle \sum\limits_{i,j=1}^n}
	a_{ij}(x_i\partial_{x_j}-y_j\partial_{y_i})
	+{\textstyle \sum\limits_{i,j=1}^n} (b_{ij}x_i\partial_{y_j}+c_{ij}y_i\partial_{x_j}) + \sigma_{\lambda-\rho}(h_{A,B,C})
\end{aligned}
\end{align}
for $A,B,C \in M_{n, n}(\mC)$ satisfying $B^{\rm T}=B$, $C^{\rm T}=C$;
\item[3)]
\begin{align}
\begin{aligned}
  \pi_\lambda(d_i)& =2\hat{z}\partial_{\hat{y}_i}+ \hat{x}_i(E_{\hat{x}}+E_{\hat{y}}+E_{\hat{z}}) +\hat{x}_i\sigma_{\lambda+\rho}(h) -{\textstyle \sum\limits_{j=1}^n} \,\hat{x}_j \sigma_{\lambda+\rho}(h_{E_{ji},0,0}) \\ &\quad  - {\textstyle \sum\limits_{j=1}^n} \,\hat{y}_j \sigma_{\lambda+\rho}(h_{0,0,E_{ij}+E_{ji}}), \\
  \pi_\lambda(e_i)& =-2\hat{z}\partial_{\hat{x}_i}+ \hat{y}_i(E_{\hat{x}}+E_{\hat{y}}+E_{\hat{z}}) +\hat{y}_i\sigma_{\lambda+\rho}(h)  + {\textstyle \sum\limits_{j=1}^n} \,\hat{y}_j \sigma_{\lambda+\rho}(h_{E_{ij},0,0}) \\ &\quad - {\textstyle \sum\limits_{j=1}^n} \,\hat{x}_j \sigma_{\lambda+\rho}(h_{0,E_{ij}+E_{ji},0}), \\
  \pi_\lambda(a)& = 4\hat{z}(E_{\hat{x}}+ E_{\hat{y}}+E_{\hat{z}}) +4\hat{z}\sigma_{\lambda+\rho}(h) - 2 {\textstyle \sum\limits_{i,j=1}^n } \hat{x}_i\hat{y}_j \sigma_{\lambda+\rho}(h_{E_{ij},0,0}) \\ &\quad + {\textstyle \sum\limits_{i,j=1}^n} \hat{x}_i\hat{x}_j \sigma_{\lambda+\rho}(h_{0,E_{ij}+E_{ji},0})  - {\textstyle \sum\limits_{i,j=1}^n} \hat{y}_i\hat{y}_j \sigma_{\lambda+\rho}(h_{0,0,E_{ij}+E_{ji}})
\end{aligned}
\end{align}
\begin{align}
\begin{aligned}
  \hat{\pi}_\lambda(d_i)& =-2y_i\partial_{z}+ \partial_{x_i}(E_x+E_y+E_z-1) -\partial_{x_i}\sigma_{\lambda-\rho}(h) +{\textstyle \sum\limits_{j=1}^n} \,\partial_{x_j} \sigma_{\lambda-\rho}(h_{E_{ji},0,0}) \\ &\quad  + {\textstyle \sum\limits_{j=1}^n} \,\partial_{y_j} \sigma_{\lambda-\rho}(h_{0,0,E_{ij}+E_{ji}}), \\
  \hat{\pi}_\lambda(e_i)& =2x_i\partial_{z}+ \partial_{y_i}(E_x+E_y+E_z-1) -\partial_{y_i}\sigma_{\lambda-\rho}(h)  - {\textstyle \sum\limits_{j=1}^n} \,\partial_{y_j} \sigma_{\lambda-\rho}(h_{E_{ij},0,0}) \\ &\quad + {\textstyle \sum\limits_{j=1}^n} \,\partial_{x_j} \sigma_{\lambda-\rho}(h_{0,E_{ij}+E_{ji},0}),  \\
  \hat{\pi}_\lambda(a)& = 4\partial_z(E_x+ E_y+E_z-1) -4\partial_z\sigma_{\lambda-\rho}(h) - 2 {\textstyle \sum\limits_{i,j=1}^n } \partial_{x_i}\partial_{y_j} \sigma_{\lambda-\rho}(h_{E_{ij},0,0}) \\ &\quad + {\textstyle \sum\limits_{i,j=1}^n} \partial_{x_i}\partial_{x_j} \sigma_{\lambda-\rho}(h_{0,E_{ij}+E_{ji},0})  - {\textstyle \sum\limits_{i,j=1}^n} \partial_{y_i}\partial_{y_j} \sigma_{\lambda-\rho}(h_{0,0,E_{ij}+E_{ji}})
\end{aligned}
\end{align}
for $i=1,2,\dots,n$.
\end{enumerate}}

\proof{The proof is a straightforward but tedious computation based on the use of commutation
relations in the Weyl algebra $\eus{A}^\mfrak{g}_{\widebar{\mfrak{u}}}$ and of the fact that
$\sigma_{\lambda+\rho} \colon \mfrak{p} \rarr \mfrak{gl}(\mathbb{V}_{\lambda+\rho})$ is a
homomorphism of Lie algebras. Let us note that a more elegant proof for a $1$-dimensional
representation of $\mfrak{p}$ easily follows from Theorem 1.3 and the proof of Theorem 3.1
in \cite{Krizka-Somberg2017}. A proof of a similar formula for an arbitrary representation
of $\mfrak{p}$ is a subject of a forthcoming work.}

Now, we shall fix a realization of the simple metaplectic submodules
of the Segal--Shale--Weil smooth admissible (unitarizable highest weight)
representation of $\mfrak{sp}(2n,\C)$ on the Schwartz
space $\mcal{S}(\R^n,\C)$. Here we restrict to its dense subspace of
$K$-finite vectors for $K={\rm U}(n)$ (the underlying Harish-Chandra module), realized
on the vector space
$\mathbb{S}=\C[\R^n] \simeq \C[q_1,q_2,\dots ,q_n]$ with the simple
metaplectic submodules of the Segal--Shale--Weil representation given
by subspaces of even and odd polynomials, respectively. The action
of generators of the metaplectic Lie algebra is
\begin{align}
\begin{aligned} \label{eq:segal-shale-weil}
  \sigma(h_{E_{ij},0,0})&=-q_j\partial_{q_i}-{\textstyle {1 \over 2}}\delta_{ij}, \\
  \sigma(h_{0,E_{ij}+E_{ji},0})&= \imag \partial_{q_i}\partial_{q_j}, \\
 \sigma(h_{0,0,E_{ij}+E_{ji}})&= \imag q_i q_j
\end{aligned}
\end{align}
for $i,j=1,2,\dots ,n$, where $\imag \in \C$ denotes the imaginary unit. The scalar product $\langle \cdot\,,\cdot \rangle \colon \mathbb{S} \otimes_\C \mathbb{S} \rarr \C$ on $\mathbb{S}$ is defined through the $\mfrak{l}^{\rm s}$-equivariant embedding into the space of Schwartz functions
$\iota \colon \mathbb{S} \rarr \mcal{S}(\R^n,\C)$, i.e.\ we have
\begin{align}
  \langle p_1,p_2 \rangle = \int_{\R^n} \iota(\widebar{p}_1)\iota(p_2)\,{\rm d}q.
\end{align}

We remark that \eqref{eq:segal-shale-weil} also defines the Segal--Shale--Weil representation on
the Schwartz space and vice versa, the action \eqref{eq:segal-shale-weil} preserves the
subspace of $K$-finite vectors in its smooth globalization on the Schwartz space.
In order to have a more efficient weight structure for our choice
of positive roots and associated Borel subalgebra, we shall work with the Fock space
realization \eqref{eq:segal-shale-weil}. This realization is also
the most convenient one from the computational point of view
and indeed, the final results are independent on the choice of a model:
by exactness of globalization functors are our results independent of the
realization we are working with.

We extend this representation of $\mfrak{l}^{\rm s}=\mfrak{sp}(2n,\C)$ to a representation
of $\mathfrak{p}$ by the trivial action of the center $\mfrak{z}(\mfrak{l})$
of $\mfrak{l}$ and by the trivial action of the nilradical $\mathfrak{u}$ of
$\mfrak{p}$, and retain the same notation $\sigma \colon \mfrak{p} \rarr \mfrak{gl}(\mathbb{S})$
for the extended action of the parabolic subalgebra
$\mfrak{p}$ of $\mfrak{g}$. In what follows, we are interested in the twisted $\mfrak{p}$-module
$\sigma_\lambda \colon \mfrak{p} \rarr \mfrak{gl}(\mathbb{S}_\lambda)$ with a twist
$\lambda \in \Hom_P(\mfrak{p},\C)$.
\medskip

Let us define the differential operators
\begin{align}
 D_s=\sum\limits_{j=1}^n \,(\imag q_j\partial_{y_j}-\partial_{x_j}\partial_{q_j}), \qquad
 E=\sum_{j=1}^{n} \,(x_{j}\partial_{x_j}+y_{j}\partial_{y_j}), \qquad
 X_s=\sum\limits_{j=1}^n \,(\imag x_jq_j+y_j\partial_{q_j}),
\end{align}
which satisfy the following commutation relations
\begin{align}
[E +n,D_s]=-D_s, \qquad [X_s,D_s]=\imag (E +n), \qquad [E +n,X_s]=X_s
\end{align}
in $\eus{A}^\mfrak{g}_{\widebar{\mfrak{u}}^*}\! \otimes_\C \End \mathbb{S}_{\lambda-\rho}$. Hence, the complex Lie algebra generated by $D_s$, $E+n$ and $X_s$ is isomorphic to $\mfrak{sl}(2,\C)$.
\medskip

\theorem{\label{thm:operator realization}
Let $\lambda \in \Hom_P(\mfrak{p},\C)$. Then the embedding of $\mfrak{g}$ into
$\eus{A}^\mfrak{g}_{\widebar{\mfrak{u}}^*}\! \otimes_\C \End \mathbb{S}_{\lambda-\rho}$ is given by
\begin{align}\label{ssw}
\begin{aligned}
  \hat{\pi}_\lambda(f_i)&=-x_i-{\textstyle {1\over 2}}z\partial_{y_i}, \\
   \hat{\pi}_\lambda(g_i)&=-y_i+{\textstyle {1\over 2}}z\partial_{x_i}, \\
  \hat{\pi}_\lambda(c)&=-z \\
\end{aligned}
\end{align}
for $i=1,2,\dots,n$;
\begin{align}\label{spn}
\begin{aligned}
 \hat{\pi}_\lambda(h)&= -E_x  - E_y - 2E_z + \lambda - (n+1), \\
   \hat{\pi}_\lambda(h_{A,B,C})&= {\textstyle \sum\limits_{i,j=1}^n}
	a_{ij}(x_i\partial_{x_j}-y_j\partial_{y_i})
	+{\textstyle \sum\limits_{i,j=1}^n} (b_{ij}x_i\partial_{y_j}+c_{ij}y_i\partial_{x_j}) + \sigma(h_{A,B,C})
\end{aligned}
\end{align}
for $A,B,C \in M_{n, n}(\C)$ satisfying $B^{\rm T}=B$, $C^{\rm T}=C$;
\begin{align}\label{soloper}
\begin{aligned}
 \hat{\pi}_\lambda(d_i)&=-2y_i\partial_z+\partial_{x_i}\big(E_x+E_y+E_z-\lambda+n-{\textstyle {1\over 2}}\big)
+q_iD_s, \\
 \hat{\pi}_\lambda(e_i)&=2x_i\partial_z+\partial_{y_i}\big(E_x+E_y+E_z-\lambda+n-{\textstyle {1 \over 2}}\big)
-\mathrm{i}\partial_{q_i}D_s,  \\
 \hat{\pi}_\lambda(a)&=4\partial_z(E_x+ E_y+E_z-\lambda +n)+ {\rm i}D_s^2
\end{aligned}
\end{align}
for $i=1,2,\dots,n$.}

\proof{The proof is a straightforward combination of Theorem \ref{reproper} and the Segal--Shale--Weil representation \eqref{eq:segal-shale-weil} twisted by a character $\lambda \in \Hom_P(\mfrak{p},\C)$.}
\vspace{-6mm}


\section{Generalized Verma modules and singular vectors}
\label{sec:singular vectors}

In what follows the generators $x_1,\dots,x_n,y_1,\dots,y_n,z$ of the graded
commutative $\C$-algebra $\C[\widebar{\mfrak{u}}^*]$ have the grading defined by
$\deg(x_i)=\deg(y_i)=1$ for $i=1,2,\dots,n$ and $\deg(z)=2$. Let us note that the choice of the grading on $(\widebar{\mfrak{u}}^*)^*$ is uniquely determined by the grading on $\widebar{\mfrak{u}}$ through the canonical isomorphism $(\widebar{\mfrak{u}}^*)^* \simeq \widebar{\mfrak{u}}$. As there is a canonical isomorphism of left $\eus{A}^\mfrak{g}_{\widebar{\mfrak{u}}^*}$-modules
\begin{align}
  \C[\widebar{\mfrak{u}}^*] \riso \eus{A}^\mfrak{g}_{\widebar{\mfrak{u}}^*}/\mcal{F}(I_e),
\end{align}
we obtain the isomorphism
\begin{align}
 \tau \circ \Phi_\lambda \colon  M^\mfrak{g}_\mfrak{p}(\mathbb{S}_{\lambda-\rho}) \riso
	\C[\widebar{\mfrak{u}}^*] \otimes_\C \mathbb{S}_{\lambda-\rho},
\end{align}
where the action of $\mfrak{g}$ on
$\C[\widebar{\mfrak{u}}^*] \otimes_\C \mathbb{S}_{\lambda-\rho}$
is given by Theorem \ref{thm:operator realization}.
Since $\C[\widebar{\mfrak{u}}^*] \otimes_\C \mathbb{S}_{\lambda-\rho}$ and
$\Sol(\mfrak{g},\mfrak{p}; \C[\widebar{\mfrak{u}}^*] \otimes_\C
\mathbb{S}_{\lambda-\rho})^\mcal{F}
\subset\C[\widebar{\mfrak{u}}^*] \otimes_\C \mathbb{S}_{\lambda-\rho}$ are semisimple
$\mfrak{l}$-modules, we denote by
$(\C[\widebar{\mfrak{u}}^*] \otimes_\C \mathbb{S}_{\lambda-\rho})_\mu$
and $\Sol(\mfrak{g},\mfrak{p};\C[\widebar{\mfrak{u}}^*]
\otimes_\C \mathbb{S}_{\lambda-\rho})^\mcal{F}_\mu$
the $\mfrak{l}$-isotypical components of highest weight $\mu \in \mfrak{h}^*$.

Let us assume $R \in \smash{\Sol(\mfrak{g},\mfrak{p};\C[\widebar{\mfrak{u}}^*]
\otimes_\C \mathbb{S}_{\lambda-\rho})^\mcal{F}_\mu}$ for some $\mu \in \mfrak{h}^*$. Then we have
\begin{align}
  \hat{\pi}_\lambda(h)R=\mu(h)R,
\end{align}
which is by Theorem \ref{thm:operator realization} equivalent to
\begin{align}
  (E_x+E_y+2E_z)R=(\lambda-(n+1)-\mu(h))R.  \label{eq:homogenity C}
\end{align}
This restricts the values of $\mu$ for which
$\Sol(\mfrak{g},\mfrak{p}; \C[\widebar{\mfrak{u}}^*]\otimes_\C
\mathbb{S}_{\lambda-\rho})_\mu^\mcal{F}$
is a non-trivial vector space, because the Euler homogeneity operators
$E_x$, $E_y$ and $E_z$ acting on
$\C[\widebar{\mfrak{u}}^*] \otimes_\C \mathbb{S}_{\lambda-\rho}$
have eigenvalues in $\N_0$. Consequently,
\begin{align}
(E_x+E_y+2E_z)R=mR \label{eq:homogenity m C}
\end{align}
for some $m \in \N_0$, and so we may apply the solution operators
\begin{align}
  \hat{\pi}_\lambda(d_i)R=0\qquad \text{and} \qquad \hat{\pi}_\lambda(e_i)R=0
\end{align}
for $i=1,2,\dots,n$ to polynomials $R$ of the form
\begin{align}
  R=\sum_{k=0}^{\lfloor {m \over 2} \rfloor} z^kR_{m-2k},
\end{align}
where $R_{m-2k} \in \C[(\mfrak{g}_{-1})^*] \otimes_\C \mathbb{S}_{\lambda-\rho}$ and satisfy
\begin{align}
(E_x + E_y)R_{m-2k}=(m-2k)R_{m-2k}. \label{eq:Rm homogenity C}
\end{align}
for $k=0,1,\dots,\lfloor {m \over 2} \rfloor$. Hence we obtain the recurrence relations
\begin{align}
  -2(k+1)y_iR_{m-2k-2}+\big(m-k-\lambda+n -{\textstyle {1\over 2}}\big)\partial_{x_i}R_{m-2k} +q_iD_sR_{m-2k}&=0, \label{eq:recurrence C1}\\
    2(k+1)x_iR_{m-2k-2}+\big(m-k-\lambda+n -{\textstyle {1\over 2}}\big)\partial_{y_i}R_{m-2k}- \imag\partial_{q_i}D_sR_{m-2k}&=0 \label{eq:recurrence C2}
\end{align}
for $k=0,1,\dots,\lfloor {m \over 2} \rfloor$, where $R_{m-2k}=0$ for $k<0$ and for $k> \lfloor {m \over 2} \rfloor$. In particular, for $k=0$ we get
\begin{align}
  -2y_iR_{m-2}+\big(m-\lambda+n -{\textstyle {1\over 2}}\big)\partial_{x_i}R_m +q_iD_sR_m&=0, \label{eq:recurrence C1 k=0}\\
  2x_iR_{m-2}+\big(m-\lambda+n -{\textstyle {1\over 2}}\big)\partial_{y_i}R_m- \imag\partial_{q_i}D_sR_m&=0. \label{eq:recurrence C2 k=0}
\end{align}

In light of the structure of the recurrence relations \eqref{eq:recurrence C1}
and \eqref{eq:recurrence C2} we see that $R$ is uniquely
determined by $R_m$. Therefore, we can define a linear mapping
\begin{align}
  \pi_{\rm top} \colon
	\Sol(\mfrak{g},\mfrak{p};\C[\widebar{\mfrak{u}}^*]\otimes_\C
	\mathbb{S}_{\lambda-\rho})^\mcal{F}_\mu
& \rarr  (\C[(\mfrak{g}_{-1})^*] \otimes_\C \mathbb{S}_{\lambda-\rho})_\mu
\end{align}
by
\begin{align}
R & \mapsto R_m,
\end{align}
which is injective and $\mfrak{l}$-equivariant.

Furthermore, the Fischer decomposition (cf.\ Appendix \ref{app:Fischer decompostion}) implies the isomorphism of vector spaces
\begin{align}
  \varphi \colon \C[\widebar{\mfrak{u}}^*] \otimes_\C \mathbb{S}_{\lambda-\rho} \riso  \bigoplus_{a,b\in \N_0} \C[z] \otimes_\C X_s^b M_a,
	\label{eq:decomposition C}
\end{align}
where $M_a = M_a^+ \oplus M_a^-$ is the subspace of $\ker D_s$ of $a$-homogeneous polynomials
in the variables $(x_1,\dots,x_n,y_1,\dots,y_n)$.
\medskip

\lemma{\label{lem:isotypic component}Let us assume $n \geq 1$. Then the isotypical components of the $\mfrak{l}$-module $\C[\widebar{\mfrak{u}}^*] \otimes_\C \mathbb{S}_{\lambda-\rho}$ are of the form
\begin{align}
  \bigoplus_{k=0}^{\lfloor \frac{a}{2}\rfloor}
	\,\C[z]_k \otimes_\C X_s^{a-2k} M_b^+ \qquad
	\text{and} \qquad \bigoplus_{k=0}^{\lfloor \frac{a}{2}\rfloor}
	\,\C[z]_k \otimes_\C X_s^{a-2k} M_b^-
\end{align}
for $a,b \in \N_0$, and they are of highest weights $\big(\lambda-(n+1)-a-2b+{1\over 2}\big)\omega_1+b\omega_2-{1\over 2}\omega_{n+1}$ and $\big(\lambda-(n+1)-a-2b+{1\over 2}\big)\omega_1+b\omega_2+\omega_n-{3\over 2}\omega_{n+1}$, respectively.}

\proof{It follows from the fact that $\C[z]_k \otimes_\C X_s^{a-2k}M_b^+$ and $\C[z]_k \otimes_\C X_s^{a-2k} M_b^-$ are the irreducible $\mfrak{l}^s$-modules with highest weights $b\omega_2-{1\over 2}\omega_{n+1}$ and $b\omega_2+\omega_n-{3\over 2}\omega_{n+1}$ for $n>1$, and $b\omega_2-{1\over 2}\omega_2$ and $b\omega_2-{3\over 2}\omega_2$ for $n=1$, respectively, cf.\ Appendix
\ref{app:Fischer decompostion} for the review of the results in \cite{DeBie-Somberg-Soucek2014}. The rest of the proof follows from \eqref{eq:homogenity C}, \eqref{eq:decomposition C} and Appendix \ref{app:Fischer decompostion}.}

\lemma{Let $m,r \in \N_0$. Then we have
\begin{align}
  D_sX_s^rv_m= -\imag{r(2m+2n+r-1) \over 2}\,X_s^{r-1}v_m. \label{DXpowers}
\end{align}
for all $v_m \in M_m$.}

\proof{By \eqref{comrel}, we have
\begin{align*}
D_sX_s^rv_m & = [D_s,X_s^r]v_m = {\textstyle \sum\limits_{j=0}^{r-1}} \,
X_s^{r-j-1}[D_s,X_s]X_s^jv_m = -\imag{\textstyle \sum\limits_{j=0}^{r-1}}\, X_s^{r-j-1}(E+n)X_s^jv_m \\
& = -\imag{\textstyle \sum\limits_{j=0}^{r-1}}\, (j+m+n)X_s^{r-1}v_m  =  -\imag{\textstyle {r(2m+2n+r-1) \over 2}}\,X_s^{r-1}v_m
\end{align*}
for all $r,m\in\mN_0$ and $v_m \in M_m$.}

\lemma{Let $m,r \in \N_0$. Then we have
\begin{align}
\partial_{x_i}X_s^rv_m &=\mathrm{i}rq_iX_s^{r-1}v_m + \mathrm{i} {r(r-1) \over 2}\, y_iX_s^{r-2}v_m
+X_s^r\partial_{x_i}v_m, \label{eq:derivative x}\\
\partial_{y_i}X_s^rv_m &=r\partial_{q_i}X_s^{r-1}v_m -\mathrm{i} {r(r-1) \over 2}\, x_iX_s^{r-2}v_m
+X_s^r\partial_{y_i}v_m  \label{eq:derivative y}
\end{align}
for all $i=1,2,\dots ,n$ and $v_m \in M_m$.}

\proof{A direct computation gives $[\partial_{x_i},X_s]={\rm i}q_i$
and $[X_s^r,q_i]=ry_iX_s^{r-1}$ for all $r\in\N_0$ and
$i=1,2,\dots ,n$. Then
\begin{align*}
\partial_{x_i}X_s^rv_m &= {\textstyle \sum\limits_{j=0}^{r-1}}\,  X_s^j[\partial_{x_i},X_s]X_s^{r-j-1} v_m +X_s^r\partial_{x_i}v_m \\
& = {\textstyle \sum\limits_{j=0}^{r-1}}\, ({\rm i}q_i X_s^{r-1}+{\rm i}[X_s^j,q_i]X_s^{r-j-1})v_m +X_s^r\partial_{x_i}v_m \\
& ={\rm i}rq_iX_s^{r-1}v_m+\mathrm{i}{\textstyle {r(r-1) \over 2}}\,y_iX_s^{r-2}v_m
+X_s^r\partial_{x_i}v_m
\end{align*}
for all $v_m \in M_m$.

The formula \eqref{eq:derivative y} follows from $[\partial_{y_i},X_s]=\partial_{q_i}$ and
$[X_s^r,\partial_{q_i}]=-\mathrm{i}rx_iX_s^{r-1}$ for $r\in\N_0$,
by analogous computation as for \eqref{eq:derivative x}. The proof is complete.}

Let us introduce the differential operators $P_1,P_2 \in \eus{A}^\mfrak{g}_{\widebar{\mfrak{u}}^*}\!\otimes_\C \End \mathbb{S}_{\lambda-\rho}$ by
\begin{align}
P_1={\textstyle \sum\limits_{j=1}^n} (x_j\hat{\pi}_\lambda(d_j) + y_j\hat{\pi}_\lambda(e_j)), \qquad
P_2={\textstyle \sum\limits_{j=1}^n} (\partial_{x_j}\hat{\pi}_\lambda(e_j)- \partial_{y_j}\hat{\pi}_\lambda(d_j)).
\end{align}

\lemma{The operators $P_1$ and $P_2$ have the following explicit form
\begin{align}
P_1 &=(E_x+E_y)\big(E_x+E_y+E_z-\lambda+n-{\textstyle {1 \over 2}}\big)-\mathrm{i}X_sD_s, \label{eq:sol operator 1 form}\\
P_2 &=2\partial_z(E_x+E_y+2n)+\mathrm{i}D_s^2. \label{eq:sol operator 2 form}
\end{align}}

\proof{It follows from \eqref{soloper}
\begin{align*}
{\textstyle \sum\limits_{j=1}^n} x_j\hat{\pi}_\lambda(d_j) &= -{\textstyle \sum\limits_{j=1}^n} 2x_jy_j\partial_z
+E_x \big(E_x+E_y+E_z-\lambda+n-{\textstyle {1 \over 2}}\big) +{\textstyle \sum\limits_{j=1}^n} x_jq_jD_s, \\
{\textstyle \sum\limits_{j=1}^n} y_j\hat{\pi}_\lambda(e_j)&= {\textstyle \sum\limits_{j=1}^n} 2x_jy_j\partial_z
+ E_y \big(E_x+E_y+E_z-\lambda+n-{\textstyle {1 \over 2}}\big)  -\mathrm{i}{\textstyle \sum\limits_{j=1}^n} y_j\partial_{q_j}D_s,
\end{align*}
and this implies \eqref{eq:sol operator 1 form}. Similarly, by \eqref{soloper} we get
\begin{align*}
-{\textstyle \sum\limits_{j=1}^n} \partial_{y_j}\hat{\pi}_\lambda(d_j) &= 2(E_y+n)\partial_z
- {\textstyle \sum\limits_{j=1}^n} \partial_{x_j}\partial_{y_j}\big(E_x+E_y+E_z-\lambda+n-{\textstyle {1 \over 2}}\big) -{\textstyle \sum\limits_{j=1}^n} \partial_{y_j}q_jD_s, \\
{\textstyle \sum\limits_{j=1}^n} \partial_{x_j}\hat{\pi}_\lambda(e_j)&= 2(E_x+n)\partial_z + {\textstyle \sum\limits_{j=1}^n} \partial_{x_j}\partial_{y_j}\big(E_x+E_y+E_z-\lambda+n-{\textstyle {1 \over 2}}\big)  -\mathrm{i}{\textstyle \sum\limits_{j=1}^n} \partial_{x_j}\partial_{q_j}D_s,
\end{align*}
and this gives \eqref{eq:sol operator 2 form}.}

Our next step is to describe the vector space of singular vectors $\Sol(\mfrak{g},\mfrak{p};\C[\widebar{\mfrak{u}}^*] \otimes_\C \mathbb{S}_{\lambda-\rho})^\mcal{F}$. Since $\smash{\Sol(\mfrak{g},\mfrak{p};\C[\widebar{\mfrak{u}}^*] \otimes_\C \mathbb{S}_{\lambda-\rho})^\mcal{F}}$ is a semisimple $\mfrak{l}$-module, first of all we decompose $\C[\widebar{\mfrak{u}}^*] \otimes_\C \!\mathbb{S}_{\lambda-\rho}$ into the isotypical components for $\mfrak{l}$, and then we find the space of solutions of $\hat{\pi}_\lambda(d_i)$ and $\hat{\pi}_\lambda(e_i)$ for $i=1,2,\dots,n$ in a given isotypical component for $\mfrak{l}$.

By Lemma \ref{lem:isotypic component} we know the form of the isotypical components of the $\mfrak{l}$-module $\C[\widebar{\mfrak{u}}^*] \otimes_\C \mathbb{S}_{\lambda-\rho}$. Let us assume $R \in \Sol(\mfrak{g},\mfrak{p};\C[\widebar{\mfrak{u}}^*] \otimes_\C \mathbb{S}_{\lambda-\rho})^\mcal{F}_\mu$ for $\mu=\big(\lambda-(n+1)-r-2m+{1\over 2}\big)\omega_1+m\omega_2-{1\over 2}\omega_{n+1}$ or $\mu= \smash{\big(\lambda-(n+1)-r-2m+{1\over 2}\big)\omega_1+m\omega_2+\omega_n-{3\over 2}\omega_{n+1}}$ and $r,m \in \N_0$, i.e.\
\begin{align}
R=\sum_{k=0}^{\lfloor {r \over 2}\rfloor}a_kz^kX_s^{r-2k}v_m, \label{singvectprop}
\end{align}
where $a_k \in \C$ for $k=0,1,\dots, \lfloor {r \over 2} \rfloor$ and either
$v_m \in M_m^+$ or $v_m \in M_m^-$. As $R \in \Sol(\mfrak{g},\mfrak{p}; \C[\widebar{\mfrak{u}}^*] \otimes_\C \mathbb{S}_{\lambda-\rho})^\mcal{F}$, we have
$R \in \ker P_1 \cap \ker P_2$ by the construction
of $P_1$ and $P_2$. Using \eqref{DXpowers}, it easily follows from
\eqref{eq:sol operator 1 form} and \eqref{eq:sol operator 2 form} that
$R \in \ker P_1 \cap \ker P_2$ if and only if
\begin{align}\label{rec1}
\big(\!(r-2k+m)(r-k+m-\lambda+n- {\textstyle {1 \over 2}})
-{\textstyle {1 \over 2}} (r-2k)(2m+2n+r-2k-1)\!
\big)a_k=0
\end{align}
and
\begin{multline}\label{rec2}
 2(k+1)(2n+m+r-2k-2)a_{k+1} \\
 -{\textstyle {\imag \over 4}} (r-2k)(r-2k-1)(2m+2n+r-2k-1)(2m+2n+r-2k-2)a_k=0
\end{multline}
hold for all $k=0,1,\dots , \lfloor {r \over 2}\rfloor$. Let us note that the coefficients $a_k \in \C$ for $k=0,1,\dots, \lfloor {r \over 2} \rfloor$ do not depend on the vector $v_m$.

The next step is to solve the recurrence relations \eqref{rec1} and \eqref{rec2}. If $a_0 = 0$, then from \eqref{rec2} we obtain $a_k = 0$ for all $k=0,1,\dots, \lfloor {r \over 2} \rfloor$.
Therefore, we may assume that $a_0 \neq 0$, and by \eqref{rec2} we have
$a_k \neq 0 $ for all $k=0,1,\dots,\lfloor {r \over 2} \rfloor$.

Now, if $m=0$ and $r=0$, there is just one equation \eqref{rec1} for $k=0$ and it is
satisfied for any $\lambda \in \C$. On the other hand, if $m \neq 0$ or $r \neq 0$, then
from \eqref{rec1} for $k=0$ we obtain
\begin{align}
  \lambda = {1 \over 2}{(m+r)^2+m(m+2n-1) \over m+r}. \label{eq:lambda}
\end{align}
The substitution for $\lambda$ from \eqref{eq:lambda} reduces the system of equations
\eqref{rec1} for $k=0,1,\dots, \lfloor {r \over 2} \rfloor$ to
\begin{align}
  {m(m+2n-1)k \over m+r}\,a_k=0 .
\end{align}
Therefore, four mutually exclusive cases have to be considered:
\begin{enumerate}
  \item[1)] $m=0$, $r=0$ and $\lambda \in \C$;
  \item[2)] $m = 0$, $r\in\mN$ and $\lambda = {1\over 2}r$;
  \item[3)] $m \neq 0$, $r=0$ and $\lambda = m+n-{1\over 2}$;
  \item[4)] $m \neq 0$, $r=1$ and $\lambda = \smash{{1 \over 2}{(m+1)^2+m(m+2n-1) \over m+1}}$.
\end{enumerate}
We shall work out each case separately, which is a content of the following lemma.
\medskip

\lemma{\label{lem:cases}
\vspace{-2mm}
\begin{enumerate}
  \item[1)] If $m=0$, $r=0$ and $\lambda \in \C$, then we have
  \begin{align}
    M_0 \subset \Sol(\mfrak{g},\mfrak{p},\C[\widebar{\mfrak{u}}^*] \otimes_\C \mathbb{S}_{\lambda-\rho})^\mcal{F}
  \end{align}
    for $\lambda \in \C$ and $n \in \N$.
  \item[2)] If $m = 0$, $r\in\mN$ and $\lambda = {1\over 2}r$, then we have
  \begin{align}
   T_{2\lambda}^n M_0 \subset \Sol(\mfrak{g},\mfrak{p},\C[\widebar{\mfrak{u}}^*]
   \otimes_\C \mathbb{S}_{\lambda-\rho})^\mcal{F}
  \end{align}
   for $\lambda = {1 \over 2} r$ and $n \in \N$, where the differential operator $T^n_a \colon \C[\widebar{\mfrak{u}}^*] \otimes_\C \mathbb{S}_{\lambda-\rho} \rarr \C[\widebar{\mfrak{u}}^*] \otimes_\C \mathbb{S}_{\lambda-\rho}$ for $a \in \N$ is defined by
  \begin{align}
   T^n_a = \sum_{k=0}^{\lfloor {a \over 2} \rfloor} \imag^k k! \binom{{a \over 2}}{k}\binom{{a \over 2}-{1\over 2}}{k}\binom{{a \over 2}-{1 \over 2}+n}{k}z^kX_s^{a-2k}. \label{eq:T operator}
  \end{align}
  \item[3)] If $m \neq 0$, $r=0$ and $\lambda = m+n-{1\over 2}$, then we have
  \begin{align}
    M_m \subset \Sol(\mfrak{g},\mfrak{p},\C[\widebar{\mfrak{u}}^*]
    \otimes_\C \mathbb{S}_{\lambda-\rho})^\mcal{F}
  \end{align}
    for $\lambda = m+n-{1\over 2}$ and $m,n \in \N$.
  \item[4)] If $m \neq 0$, $r=1$ and $\lambda = \smash{{1 \over 2}{(m+1)^2+m(m+2n-1) \over m+1}}$, then we have
  \begin{align}
    X_sM_m \subset \Sol(\mfrak{g},\mfrak{p},\C[\widebar{\mfrak{u}}^*]
    \otimes_\C \mathbb{S}_{\lambda-\rho})^\mcal{F}
  \end{align}
    for $\lambda= m+{1\over 2}$, $m\in \N$ and $n =1$.
\end{enumerate}}

\proof{1) Let us assume $m=0$, $r=0$ and $\lambda \in \C$. Then we have
\begin{align*}
  R =v_0,
\end{align*}
where either $v_0 \in M_0^+$ or $v_0 \in M_0^-$.
Since $v_0 \in \ker D_s$ and $\partial_{x_i}v_0=0$, $\partial_{y_i}v_0=0$ for $i=1,2,\dots,n$, we obtain
$\hat{\pi}_\lambda(d_i)R=0$ and $\hat{\pi}_\lambda(e_i)R=0$ for $i=1,2,\dots,n$. Therefore, we have
\begin{align*}
M_0 \subset \Sol(\mfrak{g},\mfrak{p},\C[\widebar{\mfrak{u}}^*] \otimes_\C
\mathbb{S}_{\lambda-\rho})^\mcal{F}
\end{align*}
for $\lambda \in \C$ and $n \in \N$.
\smallskip

\noindent 2) Let us assume $m=0$, $r\in\mN$ and $\lambda = {1\over 2}r$. Then we have
\begin{align*}
  R=\sum_{k=0}^{\lfloor {r \over 2}\rfloor}a_kz^kX_s^{r-2k}v_0,
\end{align*}
where either $v_0 \in M_0^+$ or $v_0 \in M_0^-$ and $a_k \in \C$ for $k=0,1,\dots,\lfloor {r \over 2} \rfloor$ satisfy the recurrence relation
\begin{align}
 2(k+1)a_{k+1} -{\textstyle {\imag \over 4}} (r-2k)(r-2k-1)(2n+r-2k-1)a_k=0. \label{eq:recurrence operator}
\end{align}
By the recurrence relations \eqref{eq:recurrence C1} and \eqref{eq:recurrence C2}, if follows that the condition $R \in \Sol(\mfrak{g},\mfrak{p}; \C[\widebar{\mfrak{u}}^*]
\otimes_\C \mathbb{S}_{\lambda-\rho})^\mcal{F}$ is equivalent to
\begin{align*}
  -2(k+1)a_{k+1}y_iX_s^{r-2k-2}v_0+\big({\textstyle {1 \over 2}r}-k+n -{\textstyle {1\over 2}}\big)a_k\partial_{x_i}X_s^{r-2k}v_0 +a_kq_iD_sX_s^{r-2k}v_0&=0, \\
  2(k+1)a_{k+1}x_iX_s^{r-2k-2}v_0+\big({\textstyle {1 \over 2}r}-k+n -{\textstyle {1\over 2}}\big)a_k\partial_{y_i}X_s^{r-2k}v_0- \imag a_k\partial_{q_i}D_sX_s^{r-2k}v_0&=0
\end{align*}
for $k=0,1,\dots, \lfloor {r \over 2} \rfloor$. Using the formulas \eqref{DXpowers}, \eqref{eq:derivative x} and \eqref{eq:derivative y}, they reduce into
\begin{align*}
  \big(\!-2(k+1)a_{k+1}+{\textstyle {\imag \over 4}}\,(r-2k)(r-2k-1)(2n+r-2k-1)a_k\big)y_iX_s^{r-2k-2}v_0&=0, \\
  \big(2(k+1)a_{k+1}-{\textstyle {\imag \over 4}}\,(r-2k)(r-2k-1)(2n+r-2k-1)a_k\big)x_iX_s^{r-2k-2}v_0&=0,
\end{align*}
which are satisfied due to the recurrence relation \eqref{eq:recurrence operator}. Let us define the differential operator $T^n_r \colon \C[\widebar{\mfrak{u}}^*] \otimes_\C \mathbb{S}_{\lambda-\rho} \rarr \C[\widebar{\mfrak{u}}^*] \otimes_\C \mathbb{S}_{\lambda-\rho}$ for $r \in \N$ by
\begin{align*}
  T^n_r = \sum_{k=0}^{\lfloor {r \over 2} \rfloor} a_k z^k X_s^{r-2k}= \sum_{k=0}^{\lfloor {r \over 2} \rfloor} \imag^k k! \binom{{r \over 2}}{k}\binom{{r \over 2}-{1\over 2}}{k}\binom{{r \over 2}-{1 \over 2}+n}{k}z^kX_s^{r-2k},
\end{align*}
where $a_k \in \C$ for $k=0,1,\dots,\lfloor {r \over 2} \rfloor$ satisfy the recurrence relation \eqref{eq:recurrence operator}. Therefore, we have
\begin{align*}
T_{2\lambda}^n M_0 \subset \Sol(\mfrak{g},\mfrak{p},\C[\widebar{\mfrak{u}}^*]
\otimes_\C \mathbb{S}_{\lambda-\rho})^\mcal{F}
\end{align*}
for $\lambda = {1 \over 2} r$ and $n \in \N$.
\smallskip

\noindent 3) Let us assume $m \neq 0$, $r=0$ and $\lambda = m+n-{1\over 2}$. Then we have
\begin{align*}
  R = v_m,
\end{align*}
where either $v_m \in M_m^+$ or $v_m \in M_m^-$. Since $v_m \in \ker D_s$ and $(E_x+E_y)v_m=mv_m$, we obtain that $\hat{\pi}_\lambda(d_i)R=0$ and $\hat{\pi}_\lambda(e_i)R=0$ for $i=1,2,\dots,n$. Therefore, we have
\begin{align*}
M_m \subset \Sol(\mfrak{g},\mfrak{p},\C[\widebar{\mfrak{u}}^*]
\otimes_\C \mathbb{S}_{\lambda-\rho})^\mcal{F}
\end{align*}
for $\lambda = m+n-{1\over 2}$ and $m,n \in \N$.
\smallskip

\noindent 4) Let us assume $m \neq 0$, $r=1$ and $\lambda = \smash{{1 \over 2}{(m+1)^2+m(m+2n-1) \over m+1}}$. Then we have
\begin{align*}
  R = X_s v_m,
\end{align*}
where either $v_m \in M_m^+$ or $v_m \in M_m^-$. From \eqref{DXpowers} we have $D_sX_sv_m=-\imag (m+n)v_m$.
By \eqref{eq:derivative x} and \eqref{eq:derivative y}, we obtain
$\partial_{x_i}X_sv_m=\imag q_i v_m + X_s\partial_{x_i}v_m$ and $\partial_{y_i}X_sv_m=\partial_{q_i} v_m + X_s\partial_{y_i}v_m$.
Hence, we may write
\begin{align*}
  \hat{\pi}_\lambda(d_i)R&=\big((m+1-\lambda+n-{\textstyle {1\over 2}})\partial_{x_i}
+q_iD_s\big)X_sv_m=\big({\textstyle {m+n \over m+1}}\,\partial_{x_i}
+q_iD_s\big)X_sv_m \\
&={\textstyle {m+n \over m+1}}\,(-\imag m q_i+X_s\partial_{x_i})v_m
\end{align*}
and
\begin{align*}
\hat{\pi}_\lambda(e_i)R&=\big((m+1-\lambda+n-{\textstyle {1\over 2}})\partial_{y_i}
-\imag \partial_{q_i}D_s\big)X_sv_m=\big({\textstyle {m+n \over m+1}}\,\partial_{y_i}
-\imag \partial_{q_i}D_s\big)X_sv_m \\
&={\textstyle {m+n \over m+1}}\,(-m \partial_{q_i}+X_s\partial_{y_i})v_m.
\end{align*}
If $n=1$, we have
\begin{align*}
  \hat{\pi}_\lambda(d_1)R&=(-\imag m q_1+X_s\partial_{x_1})v_m=(-\imag mq_1+\imag q_1x_1 \partial_{x_1} + y_1\partial_{x_1}\partial_{q_1})v_m \\
  &= (-\imag mq_1+\imag q_1x_1 \partial_{x_1} + \imag q_1 y_1\partial_{y_1})v_m =
  (-\imag mq_1+\imag q_1(E_x+E_y))v_m=0
\end{align*}
and
\begin{align*}
\hat{\pi}_\lambda(e_1)R&=(-m \partial_{q_1}+X_s\partial_{y_1})v_m = (-m \partial_{q_1}+\partial_{q_1} y_1 \partial_{y_1} + \imag x_1q_1\partial_{y_1})v_m \\
&= (-m \partial_{q_1}+\partial_{q_1} y_1 \partial_{y_1} + \partial_{q_1} x_1 \partial_{x_1})v_m = (-m \partial_{q_1}+\partial_{q_1}(E_x+E_y))v_m=0,
\end{align*}
where we used the fact that $D_sv_m=(\imag q_1 \partial_{y_1}-\partial_{x_1}\partial_{q_1})v_m=0$.

Now, if $n \neq 1$, then we may write
\begin{align*}
  \partial_{q_i}\hat{\pi}_\lambda(d_i)R&={\textstyle {m+n \over m+1}}\,(-\imag m-\imag m q_i \partial_{q_i}+\imag x_i\partial_{x_i}+X_s\partial_{x_i}\partial_{q_i})v_m, \\
   q_i\hat{\pi}_\lambda(e_i)R&={\textstyle {m+n \over m+1}}\,(-m q_i \partial_{q_i}- y_i\partial_{y_i}+X_s q_i \partial_{y_i})v_m.
\end{align*}
Putting all partial results together, we get
\begin{align*}
  {\textstyle {m+1 \over m+n}\,\sum\limits_{i=1}^n } \, ( \partial_{q_i}\hat{\pi}_\lambda(d_i)-\imag q_i\hat{\pi}_\lambda(e_i))R&=(-\imag m n+i(E_x+E_y)-X_sD_s)v_m=-\imag (n-1)mv_m \neq 0.
\end{align*}
Therefore, we have
\begin{align*}
X_sM_m \subset \Sol(\mfrak{g},\mfrak{p},\C[\widebar{\mfrak{u}}^*]
\otimes_\C \mathbb{S}_{\lambda-\rho})^\mcal{F}
\end{align*}
for $\lambda= m+{1\over 2}$, $m \in \N$ and $n =1$.}

\theorem{\label{thm:singular vector n}
\vspace{-2mm}
\begin{enumerate}
\item[1)] If $n=1$, then we have
\begin{align*}
  \tau \circ \Phi_{\lambda+\rho} \colon M^\mfrak{g}_\mfrak{p}(\mathbb{S}_\lambda)^\mfrak{u} \riso
  \begin{cases}
    M_0, & \text{if $\lambda +2 \notin {1 \over 2}\N$}, \\
    M_0 \oplus T_{2(\lambda+2)}^1M_0, & \text{if $\lambda + 1 \in \N_0$}, \\
    M_0 \oplus T_{2(\lambda+2)}^1M_0, & \text{if $\lambda=-{3 \over 2}$}, \\
    M_0 \oplus  M_{\lambda + {3 \over 2}} \oplus X_sM_{\lambda + {3 \over 2}} \oplus T_{2(\lambda+2)}^1M_0, & \text{if $\lambda + {1 \over 2} \in \N_0$}.
  \end{cases}
\end{align*}
\item[2)] If $n\geq 2$, then we have
\begin{align*}
  \tau \circ \Phi_{\lambda+\rho} \colon M^\mfrak{g}_\mfrak{p}(\mathbb{S}_\lambda)^\mfrak{u}
	\riso
  \begin{cases}
    M_0, & \text{if $\lambda + n +1 \notin {1 \over 2}\N$}, \\
    M_0 \oplus T_{2(\lambda+n+1)}^nM_0, & \text{if $\lambda + n \in \N_0$}, \\
    M_0 \oplus T_{2(\lambda+n+1)}^nM_0, & \text{if $\lambda + n + {1 \over 2} \in \N_0$, $\lambda + {1 \over 2} \notin \N_0$}, \\
    M_0 \oplus  M_{\lambda + {3 \over 2}} \oplus T_{2(\lambda+n+1)}^nM_0, & \text{if $\lambda + {1 \over 2} \in \N_0$}.
  \end{cases}
\end{align*}
\end{enumerate}}

\proof{The decomposition of the space of singular vectors $M^\mfrak{g}_\mfrak{p}(\mathbb{V})^\mfrak{u}$ is an easy consequence of Lemma \ref{lem:cases}.}

Now we illustrate the general results given in Theorem \ref{thm:singular vector n},
and write down explicit formulas for the corresponding homomorphisms between generalized Verma modules in several examples.

Let us assume $\lambda=-\big(n+1-{a \over 2}\big)\widetilde{\omega}_1$ and $\mu=-\big(n+1+{a \over 2}\big)\widetilde{\omega}_1$ for $a \in \N$. Then Theorem \ref{thm:singular vector n} implies that
\begin{align}
  T^n_a(\mathbb{S}_\lambda) \subset (\tau \circ \Phi_{\lambda+\rho})(M^\mfrak{g}_\mfrak{p}(\mathbb{S}_\lambda)^\mfrak{u})
\end{align}
for $n \in \N$, and moreover $T^n_a(\mathbb{S}_\lambda) \simeq \mathbb{S}_\mu$ as $\mfrak{p}$-modules as follows from Lemma \ref{lem:isotypic component}. Therefore, there exists a homomorphism
\begin{align}
  \varphi \colon M^\mfrak{g}_\mfrak{p}(\mathbb{S}_\mu) \rarr M^\mfrak{g}_\mfrak{p}(\mathbb{S}_\lambda)
\end{align}
of generalized Verma modules, uniquely determined by a $\mfrak{p}$-homomorphism $\varphi_0 \colon \mathbb{S}_\mu \rarr M^\mfrak{g}_\mfrak{p}(\mathbb{S}_\lambda)$ through the formula
\begin{align}
  \varphi(u\otimes v)=u\varphi_0(v)
\end{align}
for all $u \in U(\mfrak{g})$ and $v \in \mathbb{S}_\mu$. Since $\mathbb{S}_\mu$ is isomorphic
to $\mathbb{S}_\lambda$ as $\mfrak{l}^{{\rm s}}$-modules, we may set
\begin{align}
  \varphi_0(v) = (\tau\circ \Phi_{\lambda +\rho})^{-1}(T^n_av) \label{eq:phi0}
\end{align}
for all $v \in \mathbb{S}_\mu \simeq \mathbb{S}_\lambda$ (as $\mfrak{l}^{{\rm s}}$-modules) and we obtain the required homomorphism of $\mfrak{p}$-modules. Through the formula \eqref{eq:T operator} we have
\begin{align}
  T^n_{a} =  \sum_{k=0}^{\lfloor {a \over 2} \rfloor} a_k z^k X_s^{a-2k}= \sum_{k=0}^{\lfloor {a \over 2} \rfloor} \imag^k k! \binom{{a \over 2}}{k}\binom{{a \over 2}-{1\over 2}}{k}\binom{{a \over 2}-{1 \over 2}+n}{k}z^kX_s^{a-2k}. \label{eq:T operator example}
\end{align}
We denote by
\begin{align}
  Q_s = -\sum_{j=1}^n (\imag f_j q_j + g_j \partial_{q_j})
\end{align}
an element in $S(\widebar{\mfrak{u}}) \otimes_\C \End \mathbb{S}_\lambda$. Then from \eqref{eq:inverse symmetrization} and Theorem \ref{thm:operator realization} we immediately obtain
\begin{align}
  (\tau \circ \Phi_{\lambda+\rho})(\beta \otimes \id_{\mathbb{S}_\lambda})(c^m Q_s^k v)= (-1)^mz^m X_s^k v
\end{align}
for all $m,k \in \N_0$ and $v \in \mathbb{S}_\lambda$, which gives
\begin{align}
\begin{aligned}
  (\tau \circ \Phi_{\lambda+\rho})^{-1}(T^n_av) &= \sum_{k=0}^{\lfloor {a \over 2} \rfloor} (-1)^k\alpha_k (\beta\otimes \id_{\mathbb{S}_\lambda})( c^k Q_s^{a-2k}v) \\
  &= \sum_{k=0}^{\lfloor {a \over 2} \rfloor} (-1)^k\alpha_k c^k(\beta\otimes \id_{\mathbb{S}_\lambda})(Q_s^{a-2k}v) \label{eq:inverse tau Phi}
\end{aligned}
\end{align}
for all $v \in \mathbb{S}_\lambda$, where we used that $c \in \mathfrak{z}(\widebar{\mfrak{u}})$ in the last equality. Let us denote by
\begin{align}
  P_s = \sum_{j=1}^n (\imag f_j q_j + g_j \partial_{q_j})
\end{align}
an element in $U(\widebar{\mfrak{u}}) \otimes_\C \End \mathbb{S}_\lambda$. Then a straightforward computation gives
\begin{align}
  (\beta \otimes \id_{\mathbb{S}_\lambda})(Q_sv)&=P_sv, \\
  (\beta \otimes \id_{\mathbb{S}_\lambda})(Q^2_sv)&=\big(P^2_s-{\textstyle {\imag \over 2}} n c\big)v, \\
  (\beta \otimes \id_{\mathbb{S}_\lambda})(Q^3_sv)&=\big(P^3_s-{\textstyle {\imag \over 2}} (3n+1) P_s c\big)v, \\
  (\beta \otimes \id_{\mathbb{S}_\lambda})(Q^4_sv)&=\big(P^4_s-{\textstyle {\imag \over 2}} (6n+4)P_s^2c- {\textstyle {1 \over 4}}(3n^2+3n)c^2\big)v
\end{align}
for all $v \in \mathbb{S}_\lambda$. Finally, using \eqref{eq:phi0}, \eqref{eq:inverse tau Phi} and \eqref{eq:T operator example} we obtain the following explicit formulas for the homomorphisms $\varphi_0 \colon \mathbb{S}_\mu \rarr M^\mfrak{g}_\mfrak{p}(\mathbb{S}_\lambda)$ of $\mfrak{p}$-modules.

\begin{enumerate}
  \item[1)] If $a=1$, then
      \begin{align}
        \varphi_0(v)=P_s v,\quad v \in \mathbb{S}_\mu.
      \end{align}
  \item[2)] If $a=2$, then
      \begin{align}
        \varphi_0(v)=\big(P_s^2-\imag \big(n+{\textstyle {1\over 4}}\big)c\big)v, \quad v \in \mathbb{S}_\mu.
      \end{align}
  \item[3)] If $a=3$, then
      \begin{align}
        \varphi_0(v)=(P_s^3-\imag (3n+2)cP_s)v, \quad v \in \mathbb{S}_\mu.
      \end{align}
  \item[4)] If $a=4$, then
      \begin{align}
        \varphi_0(v)=\big(P_s^4-\imag \big(6n+{\textstyle {13 \over 2}}\big)cP_s^2 - \big(3n^2 + {\textstyle {9\over 2}}n + {\textstyle {9 \over 16}}\big)c^2\big)v, \quad v \in \mathbb{S}_\mu.
      \end{align}
\end{enumerate}
\vspace{-1mm}



\section{Equivariant differential operators on generalized flag manifolds}
\label{sec:invariant operators general. flag}

In the present section we retain the notation introduced in Section \ref{sec:Verma modules}.
Let us recall a few basic facts concerning infinite-dimensional complex representations of
finite-dimensional real Lie groups, see \cite{Casselman2016}.
For a complex continuous representation $(\sigma,\mV)$ of $P$, a vector $v\in\mV$ is
called smooth if the mapping $p \mapsto \sigma(p)v$ from $P$ to
$\mV$ is smooth. A continuous representation $(\sigma,\mV)$ of $P$ is smooth provided all
vectors in $\mV$ are smooth, and in particular the vector subspace $\mV^\infty\subset\mV$
of all smooth vectors is a continuous representation both for $P$ and $\mfrak{p}$.
The topological dual of a smooth representation $(\sigma,\mV)$ of $P$ is a
representation on the space of tempered distributional vectors $(\mV^*)^{-\infty}$
in the weak linear dual $\mV^*$ of $\mV$.
We define the contragredient representation of $P$ as the smooth representation $(\sigma^*,\mV^\vee)$
on the subspace of smooth vectors
$\mV^\vee=(\mV^*)^\infty \subset (\mV^*)^{-\infty}$, where the canonical
pairing $\langle\cdot\,,\cdot\rangle \colon \mV \times \mV^\vee \rarr \C$ is
non-degenerate. The previous discussion is a summary of the standard results by
Casselman, Wallach, and Schmid on the globalization of representations, cf.\ \cite{Casselman1989}, \cite{Schmid1985}.

Given a smooth $P$-module $(\sigma, \mathbb{V})$, we consider the induced representation
of $G$ on the space $\Ind_P^G(\mV)$ of smooth sections
of the homogeneous vector bundle
$\mcal{V}=G \times_P \mV \rarr G/P$ identified with
\begin{align}
{\lC}^\infty(G,\mV)^P&=\{f\in{\lC}^\infty(G,\mV);\,
f(gp)=\sigma(p^{-1})f(g)\ \text{for all}\ g\in G,\, p\in P\}.
\end{align}
We denote by $J^k_e(G,\mV)^P$ the space of $k$-jets supported at $e\in G$
of $P$-equivariant smooth mappings $f \in \mcal{C}^\infty(G,\mathbb{V})^P$ for $k\in \N_0$, and
by $J^\infty_e(G,\mV)^P$ its projective limit
\begin{align}
J^\infty_e(G,\mathbb{V})^P = {\textstyle \lim\limits_{\lrarr k}} J^k_e(G,\mV)^P.
\end{align}
A well-known fact, usually stated for a complex finite-dimensional
$P$-module $(\sigma, \mV)$, is the existence of a non-degenerate $(\mathfrak{g}, P)$-equivariant
pairing between $J^\infty_e(G,\mV)^P$  and  $M^\mfrak{g}_\mfrak{p}(\mV^*)$, see
\cite{Collingwood-Shelton1990}, \cite{Soergel1990} for independent proofs. This pairing identifies
$M^\mfrak{g}_\mfrak{p}(\mV^*)$ with the vector space of all $\C$-linear mappings
$J^\infty_e(G,\mV)^P\rarr \mC$ that factor through $J^k_e(G,\mathbb{V})^P$
for some $k \in \N_0$.
\medskip

In the following proposition we give a straightforward generalization of the statement made in the previous
paragraph for finite-dimensional inducing $P$-modules. For the reader's convenience, we give a proof of
this fact. Let us note that the proof goes along the same lines as the proof for finite-dimensional
inducing $P$-modules. However, some of the specific
results used in the proof of the claim are non-trivial for
infinite-dimensional inducing $P$-modules, and in particular, hold only for smooth inducing $P$-modules, cf.\ the manuscript \cite{Casselman2016}.

We denote by $R_g$ and $L_g$ the right and left regular action of $g \in G$ on $\mcal{C}^\infty(G,\mathbb{V})$ defined by
\begin{align}
 (R_g(f))(h)=f(hg)  \qquad \text{and} \qquad (L_g(f))(h)=f(g^{-1}h)
\end{align}
for all $h \in G$, respectively. The right and left regular representation of $G$ on $\mcal{C}^\infty(G,\mathbb{V})$ induces the right and left regular representation of the universal enveloping algebra $U(\mfrak{g})$ on $\mcal{C}^\infty(G,\mathbb{V})$ defined by
\begin{align}
  (R_X(f))(g)= {{\rm d} \over {\rm d}t}_{|t=0}f(g\exp(tX)) \qquad \text{and} \qquad (L_X(f))(g)= {{\rm d} \over {\rm d}t}_{|t=0}f(\exp(-tX)g)
\end{align}
for $X \in \mfrak{g}(\R)$, where $g \in G$, $f \in \mcal{C}^\infty(G,\mathbb{V})$, and extended to $U(\mfrak{g})$ as representations of associative $\C$-algebras, respectively. We will denote $R_u$ and $L_u$ for $u \in U(\mfrak{g})$. Moreover, we have a $\C$-bilinear mapping
\begin{align}
  \langle \cdot\,,\cdot \rangle \colon \mcal{C}^\infty(G,\mathbb{V}) \times \mathbb{V}^\vee \rarr \mcal{C}^\infty(G,\C)
\end{align}
given by
\begin{align}
  \langle f, \alpha \rangle (g) = \alpha(f(g))
\end{align}
for all $g \in G$.
\medskip

\proposition{Let $\mV$ be a smooth $P$-module. The bilinear mapping
\begin{align}
\Phi_\mV \colon \mcal{C}^\infty(G,\mV) \times (U(\mfrak{g})\otimes_\C\!\mV^\vee) \rarr \mcal{C}^\infty(G,\C)
\end{align}
defined by
\begin{align}
\Phi_\mV(f,u\otimes \alpha) = \langle R_u(f),\alpha \rangle =R_u(\langle f,\alpha \rangle),
\end{align}
where $f\in \mcal{C}^\infty(G,\mV)$, $\alpha \in\mV^\vee$ and $u\in U(\mfrak{g})$, induces a bilinear mapping
\begin{align}
\Phi_\mV \colon \mcal{C}^\infty(G,\mV)^P \times M^\mfrak{g}_\mfrak{p}(\mV^\vee) \rarr \mcal{C}^\infty(G,\C).
\end{align}
Moreover, the composition of $\Phi_\mV$ with the evaluation map at $e\in G$
\begin{align}
\begin{gathered}
  \mcal{C}^\infty(G,\mV)^P \times M^\mfrak{g}_\mfrak{p}(\mV^\vee)  \rarr \C \\
  (f,u \otimes \alpha) \mapsto \Phi_\mV(f,u \otimes \alpha)(e) \label{eq:(g,P) pairing}
\end{gathered}
\end{align}
is $(\mfrak{g},P)$-equivariant.}

\proof{Let us denote by $I(\mfrak{g}, \mfrak{p}, \mV^\vee)$ the $(\mfrak{g},P)$-submodule
of $U(\mfrak{g})\otimes_\C\! \mV^\vee$, generated by complex subspace
\begin{align*}
\langle X\otimes \alpha-1\otimes \sigma^*(X)\alpha;\,
X\in\mfrak{p},\, \alpha \in\mV^\vee \rangle \subset U(\mfrak{g})\otimes_\C\! \mV^\vee.
\end{align*}
If $f \in \mcal{C}^\infty(G,\mV)^P$ and $X \in \mfrak{p}$, then we have $R_X(f)=-\sigma(X)f$. Hence, we may write
\begin{align*}
  \Phi_\mV(f,uX\otimes \alpha)&= R_u(\langle R_X(f),\alpha \rangle)= R_u(\langle -\sigma(X)f,\alpha \rangle) \\
  &= R_u(\langle f,\sigma^*(X)\alpha \rangle)= \Phi_\mV(f,u\otimes \sigma^*(X)\alpha)
\end{align*}
for all $u \in U(\mfrak{g})$ and $\alpha \in\mV^\vee$, which implies that $\Phi_\mV(f,I(\mfrak{g},\mfrak{p},\mV^\vee))=0$ for all $f \in \mcal{C}^\infty(G,\mV)^P$.

Further, we prove that the billinear mapping \eqref{eq:(g,P) pairing} is $(\mfrak{g},P)$-equivariant. Let $f \in \mcal{C}^\infty(G,\mV)^P$, $u \in U(\mfrak{g})$ and $\alpha \in \mV^\vee$. Then we have
\begin{align*}
  \Phi_\mV(L_X(f),u \otimes \alpha)&=L_X \Phi_\mV(f,u \otimes \alpha), \\
  \Phi_\mV(f,Xu \otimes \alpha)&=R_X \Phi_\mV(f,u \otimes \alpha)
\end{align*}
for all $X \in \mfrak{g}$, where we used $[L_X,R_u]=0$. Since we have
\begin{align*}
  (R_X(f)+L_X(f))(e)=0,
\end{align*}
we get
\begin{align*}
  \Phi_\mV(L_X(f),u \otimes \alpha)(e)+ \Phi_\mV(f,Xu\otimes \alpha)(e)=0,
\end{align*}
which means that \eqref{eq:(g,P) pairing} is $\mfrak{g}$-equivariant. Furthermore, we may write
\begin{align*}
  \Phi_\mV(f,\Ad(p)u \otimes \sigma^*(p)\alpha)&=R_{\Ad(p)u}(\langle f ,\sigma^*(p)\alpha \rangle) =R_{\Ad(p)u}(\langle \sigma(p^{-1})f ,\alpha \rangle) \\
  &= R_{\Ad(p)u}(\langle R_p(f) ,\alpha \rangle) = (R_{\Ad(p)u} \circ R_p)(\langle f ,\alpha \rangle) \\
  &= (R_p \circ R_u)(\langle f ,\alpha \rangle)=R_p \Phi_\mV(f,u \otimes \alpha)
\end{align*}
for $p \in P$, where we used $R_p(f)=\sigma(p^{-1})f$ and $R_{\Ad(p)u} \circ R_p = R_p \circ R_u$. Finally, since we have
\begin{align*}
  \Phi_\mV(L_{p^{-1}}(f),u\otimes \alpha)=L_{p^{-1}}\Phi_\mV(f,u\otimes \alpha),
\end{align*}
we obtain
\begin{align*}
  \Phi_\mV(L_{p^{-1}}(f),u\otimes \alpha)(e)-\Phi_\mV(f,\Ad(p)u \otimes \sigma^*(p)\alpha)(e)=\Phi_\mV(f,u\otimes \alpha)(p)- \Phi_\mV(f,u\otimes \alpha)(p)=0,
\end{align*}
which means that \eqref{eq:(g,P) pairing} is $P$-equivariant. The proof is complete.}

Let us note that the $(\mfrak{g},P)$-equivariant bilinear mapping \eqref{eq:(g,P) pairing} identifies the generalized Verma module $M^\mfrak{g}_\mfrak{p}(\mV^\vee)$ with the vector space of all $\C$-linear mappings $J^\infty_e(G,\mV)^P \rarr \C$ that factor through $J^k_e(G,\mathbb{V})^P$ for some $k \in \N_0$. Another remark is that the duality in question is
equivalent to a more geometric construction of invariant differential jet operator.
\medskip

\corollary{\label{cor:operator-homomorphism correspondence}
Let $\mV$ and $\mW$ be smooth $P$-modules. Then there is a bijection
between $(\mfrak{g},P)$-equivariant homomorphisms of generalized Verma modules
\begin{align}
\varphi \colon  M^\mfrak{g}_\mfrak{p}(\mW^\vee) \rarr M^\mfrak{g}_\mfrak{p}(\mV^\vee),
\end{align}
and $G$-equivariant differential operators
\begin{align}
D  \colon {\lC}^\infty(G,\mV)^P \rarr {\lC}^\infty(G,\mW)^P
\end{align}
given by
\begin{align}
\langle D(f),\beta \rangle =\Phi_\mV(f, \varphi(1 \otimes \beta))
\end{align}
for $f\in \mcal{C}^\infty(G,\mV)^P$ and $\beta \in\mW^\vee$.}

We note that the proof of both injectivity and surjectivity is completely
parallel to the case of finite-dimensional inducing representations.
Now we apply this result to the case of our interest studied from the
algebraic perspective in Section \ref{sec:projective structure} and Section \ref{sec:singular vectors}.
\medskip

Let us notice that we could have defined an action of
equivariant differential operators on sections valued
in the representations given by tempered distributional
vectors rather than its dense subspace of smooth
vectors. However, the subspace of smooth vectors is
sufficient for the construction of the required duality
in Corollary \ref{cor:operator-homomorphism correspondence} and it is also more natural to work with in various
differential geometrical application. Moreover, for the
Dirac operator in symplectic geometry, the virtue of this
step was stressed in \cite{Kostant1974}.


\subsection{Equivariant differential operators on symplectic spinors in the contact projective geometry}
\label{sec:invariant operators}

By \cite{Wolf1976}, there is a double cover of the maximal parabolic subgroup
of the symplectic group $\Sp(2n+2,\mR)$ to the metaplectic group $\Mp(2n+2,\mR)$,
which splits over the unipotent radical ${\rm H}(n,\mR)$ in the
Langlands--Iwasawa decomposition of this parabolic subgroup. The group $\Mp(2n+2,\mR)$
has the maximal parabolic subgroup with the Levi subgroup isomorphic to $\GL(1,\mR)\times{\Mp(2n,\mR)}$.
We note that the extension cocycle splits over the field of complex numbers.
\medskip

Let us consider $\smash{\widetilde{G}}=\Mp(2n+2,\mR)$ for $n \in \N$, and
$\smash{\widetilde{P}}\subset\smash{\widetilde{G}}$ its maximal parabolic
subgroup with the Levi subgroup isomorphic to $\GL(1,\mR)\times{\Mp(2n,\mR)}$
and the unipotent radical in the Langlands--Iwasawa decomposition of $\smash{\widetilde{P}}$
isomorphic to the Heisenberg group ${\rm H}(n,\mR)$. As for the representation
we take the simple representations of the metaplectic group $\Mp(2n,\R)$
of the Segal--Shale--Weil representation $\mathbb{S}$, extended to
a representation of $\smash{\widetilde{P}}$ by the trivial action of $\GL(1,\R)$ and ${\rm H}(n,\R)$.
Assuming that $\lambda\in \Hom_P(\mfrak{p},\mC)$ defines a group character $e^\lambda \colon \smash{\widetilde{P}}
\rarr \GL(1,\C)$, we denote by $\mathbb{S}_\lambda$ the corresponding twisted representation
of $\smash{\widetilde{P}}$ with a twist $\lambda \in \Hom_P(\mfrak{p},\C)$. Here $\mathbb{S}_\lambda$ is understood as a smooth globalization of the corresponding Harish-Chandra module.

A non-degenerate pairing $(\cdot\,,\cdot ) \colon \mcal{S}(\R^n,\C) \otimes_\C \mathbb{S} \rarr \C$
defined by the formula
\begin{align}
  (f(q),p(q)) = \int_{\R^n} f(q)p(q)\,{\rm d}q,
\end{align}
where $q=(q_1,q_2,\dots,q_n)$, identifies $\mathbb{S}^\vee$ with $\mcal{S}(\R^n,\C)$.
The generators of $\mfrak{l}^{\rm s}$ act in the contragredient $\mfrak{p}$-module on $\mathbb{S}^\vee \simeq \mcal{S}(\R^n,\C)$
as
\begin{align}
\begin{aligned} \label{eq:segal-shale-weil dual}
  \sigma^*(h_{E_{ij},0,0})&=-q_j\partial_{q_i}-{\textstyle {1 \over 2}}\delta_{ij}, \\
  \sigma^*(h_{0,E_{ij}+E_{ji},0})&= -\imag \partial_{q_i}\partial_{q_j}, \\
 \sigma^*(h_{0,0,E_{ij}+E_{ji}})&= -\imag q_iq_j
\end{aligned}
\end{align}
for $i,j=1,2,\dots,n$, and the generators of the center $\mfrak{z}(\mfrak{l})$ of $\mfrak{l}$ and of the nilradical $\mathfrak{u}$ of $\mfrak{p}$ act trivially. Then for
$\mathbb{V}=\mathbb{S}_\lambda$ we take $\mathbb{V}^\vee=\mcal{S}(\R^n,\C)_{-\lambda}$, and denote
by $\pi^*_\lambda$ the embedding of $\mfrak{g}$ into
$\eus{A}^\mfrak{g}_{\widebar{\mfrak{u}}} \otimes_\C \End \mathbb{S}^\vee_{\lambda+\rho}$
(see Theorem \ref{reproper}) for the $\mfrak{p}$-module $(\sigma^*_{\lambda+\rho},\mathbb{S}^\vee_{\lambda+\rho})$.
\medskip

\theorem{\label{thm:edo}
The singular vectors in Theorem \ref{thm:singular vector n} (2) correspond to the $\smash{\widetilde{G}}$-equivariant
differential operators, given in the non-compact picture of the induced representations as follows.
\begin{enumerate}
\item[1)] Let $\lambda=-\big(n+1-{a \over 2}\big)\widetilde{\omega}_1$ and $\mu=-\big(n+1+{a \over 2}\big)\widetilde{\omega}_1$ for $a \in \N$.
Then there are equivariant differential operators
\begin{align}
D_a \colon {\lC}^\infty(\widebar{\mfrak{u}}(\R),\mS^\vee_{-\lambda})
\rarr  {\lC}^\infty(\widebar{\mfrak{u}}(\R),\mS^\vee_{-\mu})
\end{align}
 of order $a \in \N$. The infinitesimal intertwining property of $D_a$ is
\begin{align}
D_a\pi^*_{-{a \over 2}}(X)=\pi^*_{{a \over 2}}(X)D_a
\end{align}
for all $X \in \mfrak{g}$. We call these operators contact powers of the contact symplectic Dirac operator.
\item[2)] Let $\lambda=-\big({3 \over 2}-a\big) \widetilde{\omega}_1$ for $a\in \N$. Then there are equivariant differential operators
\begin{align}
T_a \colon {\lC}^\infty(\widebar{\mfrak{u}}(\R),\mS^\vee_{-\lambda})
\rarr  {\lC}^\infty(\widebar{\mfrak{u}}(\R),(M_{\lambda+{3 \over 2}})^\vee)
\end{align}
 of order $a \in \N$. We call these operators twistor operators on symplectic spinors.
\end{enumerate}}

\proof{The existence of $\smash{\widetilde{G}}$-equivariant differential operators easily follows from Theorem \ref{thm:singular vector n} and Corollary \ref{cor:operator-homomorphism correspondence}.}

The remaining collection of singular vectors corresponds to the \smash{$\widetilde{G}$}-equivariant differential
operator given by a multiple of the identity map. Analogous statement as the last theorem can be made in
the case $n=1$ using Theorem \ref{thm:singular vector n} (1).
\medskip

Let us define
\begin{align}
  D_{\hat{s}} ={\textstyle \sum\limits_{j=1}^n} (\imag q_j\partial_{\hat{x}_j}-\partial_{\hat{y}_j}\partial_{q_j}) \qquad \text{and} \qquad X_{\hat{s}}={\textstyle \sum\limits_{j=1}^n} (\imag \hat{y}_jq_j+\hat{x}_j\partial_{q_j}),
\end{align}
then we can write down explicit formulas for the operators $D_1, D_2, D_3, D_4$ in the form
\begin{align}
\begin{aligned}
  D_1 &= D_{\hat{s}} + {\textstyle {1 \over 2}} X_{\hat{s}} \partial_{\hat{z}}, \\
  D_2 &= \big(D_{\hat{s}} + {\textstyle {1 \over 2}} X_{\hat{s}} \partial_{\hat{z}}\big)^2 - {\textstyle {\imag \over 4}}\partial_{\hat{z}}, \\
  D_3 &= \big( D_{\hat{s}} + {\textstyle {1 \over 2}} X_{\hat{s}} \partial_{\hat{z}}\big) \big(\big(D_{\hat{s}} + {\textstyle {1 \over 2}} X_{\hat{s}} \partial_{\hat{z}}\big)^2 - \imag \partial_{\hat{z}}\big), \\
  D_4 &= \big(\big(D_{\hat{s}} + {\textstyle {1 \over 2}} X_{\hat{s}} \partial_{\hat{z}}\big)^2 - {\textstyle {\imag \over 4}}\partial_{\hat{z}}\big) \big(\big(D_{\hat{s}} + {\textstyle {1 \over 2}} X_{\hat{s}} \partial_{\hat{z}}\big)^2 - \imag {\textstyle {9 \over 4}}\partial_{\hat{z}}\big)
\end{aligned}
\end{align}
in the case $n \geq 1$.
\medskip

The examples of $\widetilde{G}$-equivariant differential operators $D_a$ for $a=1,2,3,4$
suggest the existence of the following remarkable conjectural factorization property
of $D_a$, $a \in \N$,
\begin{align}
  D_{2k+2}&= {\textstyle \prod\limits_{j=0}^k} \big(\big(D_{\hat{s}} + {\textstyle {1 \over 2}} X_{\hat{s}} \partial_{\hat{z}}\big)^2-\imag {\textstyle {(2j+1)^2 \over 4}}\,\partial_{\hat{z}}\big), \\
  D_{2k+1}&=\big(D_{\hat{s}} + {\textstyle {1 \over 2}} X_{\hat{s}} \partial_{\hat{z}}\big) {\textstyle \prod\limits_{j=1}^k} \big(\big(D_{\hat{s}} + {\textstyle {1 \over 2}} X_{\hat{s}} \partial_{\hat{z}}\big)^2-\imag j^2 \partial_{\hat{z}}\big)
\end{align}
for $k \in \N_0$. We do not know a direct proof of this observation.


\begin{appendices}

\section{The Fischer decomposition for $\mfrak{sp}(2n,\C)$}
\label{app:Fischer decompostion}

Throughout the article, we needed a rather precise information about the
simple part $\mfrak{l}^s\simeq\mathfrak{sp}(2n,\mC)$ of the
Levi subalgebra $\mathfrak{l}$-module structure on the generalized
Verma module $M^\gog_\gop(\mS_{\lambda-\rho})$ for $\lambda\in\Hom_P(\mfrak{p},\C)$.
Let us first mention that as discussed beyond Definition \ref{vmgendef} in
Section \ref{sec:Verma modules}, the present decomposition holds
for both the Harish-Chandra module $\mS_{\lambda-\rho}$ and its smooth globalization.
It can be described in terms of the Howe duality for the pair
$\mathfrak{sl}(2,\mC) \oplus \mathfrak{sp}(2n,\mC)$ acting on
$\C[\widebar{\mfrak{u}}^*] \otimes_\C \mathbb{S}_{\lambda-\rho}$ by \eqref{spn} for $\mfrak{sp}(2n,\C)$ and by
\begin{align}
 D_s=\sum\limits_{j=1}^n \,(\imag q_j\partial_{y_j}-\partial_{x_j}\partial_{q_j}), \qquad
 E=\sum_{j=1}^{n} \,(x_{j}\partial_{x_j}+y_{j}\partial_{y_j}), \qquad
 X_s=\sum\limits_{j=1}^n \,(\imag x_jq_j+y_j\partial_{q_j})
\end{align}
for $\mfrak{sl}(2,\C)$. For a detailed exposition, see \cite{DeBie-Somberg-Soucek2014}.
This decomposition can be schematically represented as
\begin{align}
\C[(\mfrak{g}_{-1})^*]\otimes \mS_{\lambda-\rho}
\simeq \smash{\bigoplus_{a,b \in \N_0}} X_s^a M_b,
\quad M_b=\big(\C[(\mfrak{g}_{-1})^*]_b\otimes \mS_{\lambda-\rho}\big)
\cap\ker D_s.
\end{align}
\begin{eqnarray}\label{obrdekonposition}\nonumber
\xymatrix@=11pt{P_0 \otimes \mathbb{S} \ar@{=}[d] &  P_1 \otimes \mathbb{S} \ar@{=}[d]&
P_2 \otimes \mathbb{S} \ar@{=}[d] & P_3 \otimes \mathbb{S} \ar@{=}[d] &
P_4 \otimes \mathbb{S} \ar@{=}[d]& P_5 \otimes \mathbb{S} \ar@{=}[d] \\
M_0 \ar[r] & X_s M_0 \ar @{} [d] |{\oplus} \ar[r] & X_s^2 M_0 \ar @{} [d] |{\oplus} \ar[r] & X_s^3 M_0 \ar @{} [d] |{\oplus}
 \ar[r] & X_s^4 M_0 \ar @{} [d] |{\oplus}\ar[r] & X_s^5 M_0 \ar @{} [d] |{\oplus} \\
& M_1 \ar[r] & X_s M_1 \ar @{} [d] |{\oplus}\ar[r] & X_s^2 M_1 \ar @{} [d] |{\oplus}
 \ar[r] & X_s^3 M_1 \ar @{} [d] |{\oplus}\ar[r] & X_s^4 M_1 \ar @{} [d] |{\oplus} \\
&& M_2 \ar[r] & X_s M_2 \ar @{} [d] |{\oplus}
 \ar[r] & X_s^2 M_2 \ar @{} [d] |{\oplus}\ar[r] & X_s^3 M_2 \ar @{} [d] |{\oplus} \\
&&& M_3 \ar[r] & X_s M_3 \ar @{} [d] |{\oplus}\ar[r] & X_s^2 M_3  \ar @{} [d] |{\oplus} \\
&&&& M_4 \ar[r] & X_s M_4 \ar @{} [d] |{\oplus}  \\
&&&&& M_5}
\end{eqnarray}
We used the notation $P_b$ for the $b$-homogeneous polynomials
instead of $\C[(\mfrak{g}_{-1})^*]_b$ in the last picture.
The operators $D_s$ and $X_s$ act in
the previous picture horizontally, but in the opposite
direction, $E$ preserves each simple symplectic
module in the decomposition and the $\mathfrak{sl}(2,\mC)$-commutation relations
in $\eus{A}^\mfrak{g}_{\widebar{\mfrak{u}}^*}\!\otimes_\C \End \mathbb{S}_{\lambda-\rho}$ are
\begin{align}\label{comrel}
[E +n,D_s]=-D_s, \qquad [X_s,D_s]=\imag (E +n), \qquad [E +n,X_s]=X_s.
\end{align}
For the operator $a \otimes \id_{\mathbb{S}_{\lambda-\rho}} \in \eus{A}^\mfrak{g}_{\widebar{\mfrak{u}}^*}\!\otimes_\C \End \mathbb{S}_{\lambda-\rho}$
with $a \in \C$, we use the shorthand notation $a$. Let us note that the operators $X_s$, $D_s$ and $E$ commute with the action of $\mathfrak{sp}(2n,\mC)$ on $\C[\widebar{\mfrak{u}}^*] \otimes_\C \mathbb{S}_{\lambda-\rho}$.

\end{appendices}


\section*{Acknowledgement}

L.\,Křižka is supported by PRVOUK p47,
P.\,Somberg acknowledges the financial support from the grant GA\,CR P201/12/G028.



\begin{thebibliography}{10}

\bibitem{Bernstein-Gelfand1980}
Joseph~N. Bernstein and Sergei~I. Gelfand, \emph{{Tensor products of finite-
  and infinite-dimensional representations of semisimple Lie algebras}},
  Compositio Math. \textbf{41} (1980), no.~2, 245--285.

\bibitem{Bernstein-Krotz2014}
Joseph~N. Bernstein and Bernhard Krötz, \emph{{Smooth Fréchet globalizations of
  Harish-Chandra modules}}, Israel J. Math. \textbf{199} (2014), no.~1,
  45--111.

\bibitem{DeBie-Somberg-Soucek2014}
Hendrik~De Bie, Petr Somberg, and Vladimír Souček, \emph{{The metaplectic Howe
  duality and polynomial solutions for the symplectic Dirac operator}}, J.
  Geom. Phys. \textbf{75} (2014), 120--128.

\bibitem{Casselman1989}
William~A. Casselman, \emph{{Canonical extensions of Harish-Chandra modules to
  representations of $G$}}, Canad. J. Math. \textbf{41} (1989), no.~3,
  385--438.

\bibitem{Casselman2016}
Bill Casselman, \emph{{Continuous representations}}, available at {\tt
  https://www.math.ubc.ca/$\sim$cass/\\research/pdf/Continuous.pdf}, 2016.

\bibitem{Collingwood-Shelton1990}
David~H. Collingwood and Brad Shelton, \emph{{A duality theorem for extensions
  of induced highest weight modules}}, Pacific J. Math. \textbf{146} (1990),
  no.~2, 227--237.

\bibitem{Fox2005}
Daniel J.~F. Fox, \emph{{Contact projective structures}}, Indiana Univ. Math.
  J. \textbf{54} (2005), no.~6, 1547--1598.

\bibitem{Gover-Graham2005}
Rod~A. Gover and Robin~C. Graham, \emph{{CR invariant powers of the
  sub-Laplacian}}, J. Reine Angew. Math. \textbf{583} (2005), 1--27.

\bibitem{Graham-Jenne-Mason-Sparling1992}
Robin~C. Graham, Ralph Jenne, Lionel~J. Mason, and George~A. Sparling,
  \emph{{Conformally invariant powers of the Laplacian. I. Existence}}, J.
  London Math. Soc. (2) \textbf{46} (1992), no.~3, 557--565.

\bibitem{Habermann-Habermann2006}
Katharina Habermann and Lutz Habermann, \emph{{Introduction to symplectic Dirac
  operators}}, Lecture Notes in Mathematics, vol. 1887, Springer-Verlag,
  Berlin, 2006.

\bibitem{Holland-Sparling2001}
Jonathan Holland and George Sparling, \emph{{Conformally invariant powers of
  the ambient Dirac operator}}, {\tt arXiv:math/0112033} (2001).

\bibitem{Kadlcakova2001}
Lenka Kadlčáková, \emph{{Dirac operators in parabolic contact symplectic
  geometry}}, Ph.D. thesis, Charles University in Prague, 2001.

\bibitem{Krizka-Somberg2017}
Libor Křižka and Petr Somberg, \emph{{Algebraic analysis of scalar generalized
  Verma modules of Heisenberg parabolic type I: $A_n$-series}}, Transform.
  Groups (2017), {\tt doi:10.1007/s00031-016\\-9414-5}.

\bibitem{Krizka-Somberg2016}
\bysame, \emph{{Algebraic analysis on scalar generalized Verma modules of
  Heisenberg parabolic type II.: $C_n, D_n$-series}},  (in preparation).

\bibitem{koss}
Toshiyuki Kobayashi, Bent \O{}rsted, Petr Somberg, and Vladimír Souček,
  \emph{{Branching laws for Verma modules and applications in parabolic
  geometry. I}}, Adv. Math. \textbf{285} (2015), 1796--1852.

\bibitem{Kobayashi-Pevzner2016a}
Toshiyuki Kobayashi and Michael Pevzner, \emph{{Differential symmetry breaking
  operators: I. General theory and F-method}}, Selecta Math. (N.S.) \textbf{22}
  (2016), no.~2, 801--845.

\bibitem{Kostant1974}
Bertram Kostant, \emph{{Symplectic spinors}}, Symposia Mathematica, Progress in
  Mathematics, vol. XIV, Academic Press, London, 1974, pp.~139--152.

\bibitem{Kostant1975}
\bysame, \emph{{Verma modules and the existence of quasi-invariant differential
  operators}}, {Non-Commutative Harmonic Analysis}, Lecture Notes in
  Mathematics, vol. 466, Springer, Berlin, 1975, pp.~101--128.

\bibitem{Kriegl-Michor1997}
Andreas Kriegl and Peter~W. Michor, \emph{{The convenient setting of global
  analysis}}, Mathematical Surveys and Monographs, vol.~53, American
  Mathematical Society, Providence, 1997.

\bibitem{Mazorchuk-Stroppel2008}
Volodymyr Mazorchuk and Catharina Stroppel, \emph{{Categorification of
  (induced) cell modules and the rough structure of generalised Verma
  modules}}, Adv. Math. \textbf{219} (2008), no.~4, 1363--1426.

\bibitem{Schmid1985}
Wilfried Schmid, \emph{{Boundary value problems for group invarint differential
  equations}}, The mathematical heritage of Élie Cartan (Lyon, 1984),
  Astérisque, Numéro Hors Série, Société Mathématique de France, Paris,
  1985, pp.~311--321.

\bibitem{Soergel1990}
Wolfgang Soergel, \emph{{The prime spectrum of the enveloping algebra of a
  reductive Lie algebra}}, Math. Z. \textbf{204} (1990), no.~4, 559--581.

\bibitem{Wolf1976}
Joseph~A. Wolf, \emph{{Unitary representations of maximal parabolic subgroups
  of the classical groups}}, vol.~8, Mem. Amer. Math. Soc., no. 180, American
  Mathematical Society, Providence, 1976.

\end{thebibliography}

\providecommand{\bysame}{\leavevmode\hbox to3em{\hrulefill}\thinspace}
\providecommand{\MR}{\relax\ifhmode\unskip\space\fi MR }
\providecommand{\MRhref}[2]{%
  \href{http://www.ams.org/mathscinet-getitem?mr=#1}{#2}
}
\providecommand{\href}[2]{#2}

\end{document}